\documentclass[12pt]{article}

\usepackage[margin=1in]{geometry}
\usepackage{xcolor}

\usepackage{amsmath,amsthm,amssymb,mathtools,bbm,mathrsfs}
\usepackage[authoryear,longnamesfirst,square]{natbib}
\renewcommand{\cite}{\citep}

\usepackage{array}
\usepackage{enumitem}
\usepackage{textcase}

\usepackage[breaklinks=true]{hyperref}
\hypersetup{
colorlinks=true,
linkcolor=black,
urlcolor=blue,
citecolor=black
}

\numberwithin{equation}{section}
\allowdisplaybreaks[4]

\newtheoremstyle{plain2}{\topsep}{\topsep}{\itshape}
{0pt}{\bfseries}{.}{.5em}{}
\newtheoremstyle{definition2}{\topsep}{\topsep}{}
{0pt}{\bfseries}{.}{.5em}{}

\theoremstyle{plain2}
\newtheorem{theorem}{Theorem}[section]
\newtheorem{proposition}[theorem]{Proposition}
\newtheorem{lemma}[theorem]{Lemma}

\theoremstyle{definition2}
\newtheorem{definition}[theorem]{Definition}

\newtheorem{remark}[theorem]{Remark}

\makeatletter

\renewcommand{\phi}{\varphi}
\newcommand{\eps}{\varepsilon}
\newcommand{\eq}{\eqref}

\newcommand{\dlp}{\mathop{d_{\mathrm{LP}}}}

\newcommand{\bigo}{\mathrm{O}}

\def\tsfrac#1#2{{\textstyle\frac{#1}{#2}}}

\newcommand{\Po}{\mathop{\mathrm{Po}}}

\newcommand{\law}{\mathscr{L}}
\newcommand{\eqlaw}{\stackrel{d}{=}}

\newcounter{ctr}\loop\stepcounter{ctr}\edef\X{\@Alph\c@ctr}%
	\expandafter\edef\csname s\X\endcsname{\noexpand\mathscr{\X}}
	\expandafter\edef\csname c\X\endcsname{\noexpand\mathcal{\X}}
	\expandafter\edef\csname b\X\endcsname{\noexpand\boldsymbol{\X}}
	\expandafter\edef\csname I\X\endcsname{\noexpand\mathbb{\X}}
	\expandafter\edef\csname r\X\endcsname{\noexpand\mathrm{\X}}
\ifnum\thectr<26\repeat

\def\be#1{\begin{equation*}#1\end{equation*}}
\def\ben#1{\begin{equation}#1\end{equation}}
\def\bes#1{\begin{equation*}\begin{split}#1\end{split}\end{equation*}}
\def\besn#1{\begin{equation}\begin{split}#1\end{split}\end{equation}}

\def\ba#1{\begin{align*}#1\end{align*}}
\def\ban#1{\begin{align}#1\end{align}}
\def\given{\mskip 0.5mu plus 0.25mu\vert\mskip 0.5mu plus 0.15mu}
\newcounter{@bracketlevel}
\def\@bracketfactory#1#2#3#4#5#6{
\expandafter\def\csname#1\endcsname##1{%
\addtocounter{@bracketlevel}{1}%
\global\expandafter\let\csname @middummy\alph{@bracketlevel}\endcsname\given%
\global\def\given{\mskip#5\csname#4\endcsname\vert\mskip#6}\csname#4l\endcsname#2##1\csname#4r\endcsname#3%
\global\expandafter\let\expandafter\given\csname @middummy\alph{@bracketlevel}\endcsname
\addtocounter{@bracketlevel}{-1}}%
}
\def\bracketfactory#1#2#3{%
\@bracketfactory{#1}{#2}{#3}{relax}{0.5mu plus 0.25mu}{0.5mu plus 0.15mu}
\@bracketfactory{b#1}{#2}{#3}{big}{1mu plus 0.25mu minus 0.25mu}{0.6mu plus 0.15mu minus 0.15mu}
\@bracketfactory{bb#1}{#2}{#3}{Big}{2.4mu plus 0.8mu minus 0.8mu}{1.8mu plus 0.6mu minus 0.6mu}
\@bracketfactory{bbb#1}{#2}{#3}{bigg}{3.2mu plus 1mu minus 1mu}{2.4mu plus 0.75mu minus 0.75mu}
\@bracketfactory{bbbb#1}{#2}{#3}{Bigg}{4mu plus 1mu minus 1mu}{3mu plus 0.75mu minus 0.75mu}
}
\bracketfactory{clc}{\lbrace}{\rbrace}
\bracketfactory{clr}{(}{)}
\bracketfactory{cls}{[}{]}
\bracketfactory{abs}{\lvert}{\rvert}
\bracketfactory{norm}{\Vert}{\Vert}
\bracketfactory{floor}{\lfloor}{\rfloor}
\bracketfactory{ceil}{\lceil}{\rceil}
\bracketfactory{angle}{\langle}{\rangle}

\def\leq{\mathchoice%
{\mathrel{\mkern6mu\mathchar"0436\mkern6mu}}%
{\mathchar"3436}{\mathchar"3436}{\mathchar"3436}}
\let\le\leq
\def\geq{\mathchoice%
{\mathrel{\mkern6mu\mathchar"043E\mkern6mu}}%
{\mathchar"343E}{\mathchar"343E}{\mathchar"343E}}
\let\ge\geq

\newcount\minute
\newcount\hour
\newcount\hourMins
\def\now{%
\minute=\time%
\hour=\time \divide \hour by 60%
\hourMins=\hour \multiply\hourMins by 60%
\advance\minute by -\hourMins%
\zeroPadTwo{\the\hour}:\zeroPadTwo{\the\minute}%
}
\def\zeroPadTwo#1{\ifnum #1<10 0\fi#1}

\renewcommand\section{\@startsection {section}{1}{\z@}%
{-3.5ex \@plus -1ex \@minus -.2ex}%
{1.3ex \@plus.2ex}
{\center\small\sc\mathversion{bold}\MakeTextUppercase}}

\def\subsection#1{\@startsection {subsection}{2}{0pt}%
{-3.5ex \@plus -1ex \@minus -.2ex}%
{1ex \@plus.2ex}%
{\bf\mathversion{bold}}{#1}}

\def\subsubsection#1{\@startsection{subsubsection}{3}{0pt}%
{\medskipamount}%
{-10pt}%
{\normalsize\itshape}{\kern-2.2ex. #1.}}

\def\blfootnote{\xdef\@thefnmark{}\@footnotetext}

\makeatother

\def\sp#1{^{(#1)}}
\def\st#1{^{[#1]}}

\def\a{\alpha}
\def\b{\beta}

\def\s{\sigma}
\def\d{\delta}

\def\k{\kappa}

\def\m{\mu}

\def\law{{\mathcal L}}

\def\ignore#1{}

\def\tU{{\widetilde U}}
\def\hZ{{\widehat Z}}

\def\hJ{{\widehat J}}
\def\hX{{\widehat X}}
\def\hC{{\widehat C}}

\def\uii{^{(i)}}
\def\ui{^{(1)}}
\def\ut{^{(2)}}
\def\uts{^{[2]}}
\def\uh{^{(3)}}

\def\uo{^{(0)}}

\def\h{\eta}

\def\tG{{\widetilde G}}

\def\l{\lambda}
\def\th{\theta}

\def\ee{\epsilon}
\def\b{\beta}

\def\f{\varphi}
\def\ff{\phi}
\def\Eq{\ =\ }
\def\Def{\ :=\ }
\def\Ref{\eqref}
\def\ex{\IE}
\def\pr{\IP}

\newcommand{\E}{\mathbb{E}}

\renewcommand{\le}{\leqslant}
\renewcommand{\leq}{\leqslant}
\renewcommand{\ge}{\geqslant}
\renewcommand{\geq}{\geqslant}

\newcommand{\N}{\mathbb{N}}

\newcommand{\R}{\mathbb{R}}

\def\w{w}
\def\wt{\widetilde}
\def\tX{{\wt X}}
\def\tN{{\wt N}}

\def\tY{{\wt Y}}
\def\tC{{\wt C}}
\def\wT{{\widehat T}}
\def\hM{{\widehat M}}
\def\hMM{M^*}
\def\hNN{N^*}
\def\hXX{X^*}

\def\thXX{\tX^*}
\def\Le{\ \leq\ }
\def\z{\zeta}
\def\MM{M}
\def\non{\nonumber}
\def\hM{{\widehat M}}
\def\cM{{\mathcal M}}
\def\hT{\mathrm{T}}
\def\bone{{\mathbf 1}}
\def\tbe{{\b}}
\def\bo{{\b}}

\def\LL{\h}
\def\oG{{\overline G}}
\def\bab{{\bar\b}}
\def\stin{}
\def\tr{\top}
\def\boe{\mathbf{e}}
\def\tg{{\tilde g}}
\def\ppsi{\Psi}

\begin{document}

\title{\sc\bf\large\MakeUppercase{Stein's method, Gaussian processes and Palm measures, with applications to queueing}}
\author{\sc A.\ D.\ Barbour, Nathan Ross, Guangqu Zheng}
\date{\it Universit\"at Z\"urich, University~of~Melbourne, University of Edinburgh}
\maketitle

\begin{abstract} 
We develop a general approach to Stein's method for approximating a random process in the path space 
$\ID([0,\hT]\to \IR^d)$ by a real continuous  Gaussian 
process.  We then use the approach in the context of processes that have a representation as integrals with respect to an
underlying point process, deriving a general \emph{quantitative} Gaussian approximation.  
The error bound is expressed  in terms of couplings of the original process to processes generated from the reduced 
Palm measures associated with the  point process.  As applications, we study  
certain $\text{GI}/\text{GI}/\infty$ queues in the ``heavy traffic'' regime. 
\end{abstract}


\section{Introduction}

Gaussian processes arise as approximations to real processes in a wide variety of applications. Often, the approximation is 
taken as read, and Gaussian processes become part of the model, as in stochastic integrals in finance.  In other circumstances, 
as in queuing systems, they arise as approximations in the limit; see, for example, \cite{Robert2003} and \cite{Pang2007}. 
The fundamental example, which forms the basis of many other limiting results,
is Donsker's theorem, which states that random walk, after proper normalization, converges weakly in path space to  
Brownian motion. 
Then probabilities for systems that converge to a Gaussian process may be approximated by the analogous limiting probabilities, 
which are typically more tractable, due to the many beautiful and useful properties of Gaussians.  A key task in this 
setting is to estimate the error made in the approximation. For Donsker's theorem, this is well understood, but for 
more general processes there are few results. 

In this paper, we 
establish a Stein equation, together with properties of its solutions, suitable for use in quantifying the error
in approximating a multi-dimensional c\`adl\`ag process by a general Gaussian process. 
The approach generalizes and improves the theory for approximation
by Brownian motion presented in \cite{Barbour1990}, and dovetails with \cite{Barbour2021} to give bounds on the error in 
terms of the L\'evy--Prokhorov distance, which metrizes weak convergence with respect to the Skorokhod topology.  
As a concrete application of the method,
we prove a general result, Theorem~\ref{thm:stnpalm}, that gives such bounds when the process being approximated 
can be expressed as an integral with respect to a point process.  Theorem~\ref{thm:stnpalm} is then applied to  
$\text{M}/\text{GI}/\infty$ and $\text{GI}/\text{GI}/\infty$ queues in the 
heavy traffic regime, obtaining the first rates of convergence in some settings that are closely related 
to limiting approximations given in
\cite{Iglehart1965, Borovkov1967, Whitt1982, Krichagina1997}
and \cite{Puhalskii2010},  where the limiting processes are typically not Brownian motion.

\subsection{Setup}\label{sec:setup}
Let $N$ be a  simple  point process on $\IR^d$ with mean measure
$\lambda$. 
 Let the collection of functions 
\be{
   \{J_u\colon\,[0,\hT]\to \IR^{ {p}}\}_{u\in \IR^d},
}
be such that $J_{\cdot}(s) \in L^2(\IR^d {\to\IR^p}, \lambda)$ for all $s\in[0,\hT]$, and $(u,s)\longmapsto J_{u}(s)$ is jointly 
measurable. 
We study the random process
$X:[0,\hT]\to \mathbb{R}^{ {p}}$ defined by
\be{
    X(s) \Def \int_{\IR^d} J_{u}(s) N(du).  
}
 {It is shown in Sections~\ref{sec:mginfres} and~\ref{sec:mggres} that  many queueing processes can be written in this form. 
Before going into specific detail, we first need to establish the general framework.} 

Define the centered and scaled random measure $\tN:= \sigma^{-1}(N-\lambda)$, where  $\sigma>0$ is a scaling 
parameter, and define the process $\tX$ by
\ben{\label{eq:gentildex}
   \tX(s) \Def \int_{\IR^d}  J_{u}(s) \tN(du).
}
We are interested in the distribution of $\tX$ when $\lambda$ is  at  ``high intensity'', meaning that the mass of $\lambda$ 
is large, and the choice of $\sigma$  stabilizes $\widetilde{N}$; and, in particular, we want to approximate the distribution 
of $\widetilde X$ by that of a Gaussian process $(Z(t),\, t\in[0,\hT])$.
Informally, we think of  $\tN$ as close in distribution to a centered Gaussian  random measure~$\rN$ with 
intensity measure $\Lambda \approx \sigma^{-2}\lambda$, and then in turn $Z(s)=\int_{\IR^d}\hJ_u(s) \rN(du)$, 
for some possibly different family of functions $\{\hJ_u\}_{u\in\IR^d}$. Formally, $Z$ is  a centered Gaussian 
process with covariance function
\be{
    \IE[Z(s) Z(t) {^{\tr}}] \Eq \int_{\IR^d} \hJ_u(s) \hJ_u(t)^{ {\tr}} \Lambda(du),\,\,  t,s\in[0,\hT],
}
 {where {$Z$ and} $\hJ$ are column vectors, and  $\tr$ denotes transpose.}
If $N$ is a Poisson process, then 
the natural approximating Gaussian process has $\Lambda =\sigma^{-2} \lambda$ and $\hJ_u=J_u$, but this is not 
necessarily the case for other point processes.
 
The corresponding approximation result, Theorem~\ref{thm:stnpalm} below, gives a bound on 
\be{
              \babs{\IE\bcls{ g(\tX)}-\IE\bcls{ g(Z)}} 
}
for a certain set~$\cM$ of ``test" functions~$g$. 
The bound in the approximation result is completely general
for~$X$ of the form above, but requires the construction of close couplings $(N,N\sp{u})_{u\in \IR^d}$ of $N$ with 
its ``reduced Palm measures'' at $u\in\IR^d$. In the case where~$N$ is a Poisson process, we can set $N\sp{u}=N$, 
and our bound becomes very simple;  see \eqref{eq:poisstnbd}.  The test functions are described 
in detail in Section~\ref{sec:testfunc}, but 
are essentially those introduced in \cite{Barbour1990}, and include smooth functions of the process at a fixed 
number of times. Such test functions are now commonly used for Gaussian process approximation in the Stein's 
method literature; see, for example \cite{Kasprzak2017a, Kasprzak2020a}. For a sequence of processes 
$(\tX_n)_{n\geq1}$,  {the fact that} $\abs{\IE\cls{ g(\tX_n)}-\IE\cls{ g(Z)}} \to 0$ for all test functions~$g$  in~$\cM$
does not alone imply weak convergence of the processes with respect to either the supremum or the Skorokhod topologies, 
but with the results of 
 \cite{Barbour2021} and a little extra work, it is not too difficult in our applications to obtain bounds on 
the L\'evy--Prokhorov distance (with respect to the Skorokhod topology) that tend to zero, and hence imply 
weak convergence.

We next discuss the test functions in detail.

\subsection{Test functions}\label{sec:testfunc}

Let  $\ID^{ {p}}:=\ID {\bclr{[0,\hT] \to \IR^p}}$ be the set of functions from $[0,\hT]$ to $\IR^{ {p}}$ that are right continuous with left limits. 
Endowed with the sup norm, $\ID {^p}$ is a Banach space (though not separable), and so for functions
 $g\colon\,\ID {^p}\to\IR$, we denote their Fr\'echet derivatives  by~$D, D^2,\ldots$. 
Following \cite{Barbour1990} (see also \cite{Kasprzak2017}), for $g\colon\,\ID^{ {p}}\to\IR$,  we define
\be{
  \norm{g}_L \Def \sup_{w\in \ID {^p}}\frac{\abs{g(w)}}{1+\norm{w}^3},
}
where $\norm{w}=\sup_{0\leq t\leq \hT} \abs{w(t)}$ denotes the sup-norm,
and then define the Banach space
\be{
   L \Def \bbclc{g\colon\,\ID^{ {p}}\to \IR \colon\ \mbox{$g$ is continuous and } \norm{g}_L<\infty}.
}
For $g$ twice Fr\'echet differentiable, we define 
\be{
   \norm{g}_M \Def \norm{g}_L+\sup_{w\in \ID^{ {p}}}\frac{\norm{D g(w)}}{1+\norm{w}^2}
           +\sup_{w\in \ID^{ {p}}}\frac{\norm{D^2 g(w)}}{1+\norm{w}}
                 +\sup_{w,h\in \ID^{ {p}}}\frac{\norm{D^2 g(w+h)-D^2 g(w)}}{\norm{h}}\,,
}
where $\norm{A}:=\sup_{w:\norm{w}=1} \abs{A[ {w^{[k]}]}}$ for~$A$ a $k$-linear form, and 
$A[ {w^{[k]}]} := A[w,w,\ldots,w]$. This leads to the space
\be{
    M \Def \bbclc{g\colon\,\ID^{ {p}}\to \IR \colon\ \mbox{$g$ is twice Fr\'echet differentiable and } \norm{g}_M<\infty}.
}
We also work on its subspace 
\be{
    M'\Def \bbclc{g\in M\colon\ \sup_{w\in \ID^{ {p}}} \norm{D^2 g(w)}< \infty },
}
and for $g\in M'$, we  define the norm
\be{
\norm{g}_{M'} \Def
 \sup_{w\in \ID^{ {p}}}\frac{\abs{g(w)}}{1+\norm{w}^2}
 +\sup_{w\in \ID^{ {p}}}\frac{\norm{D g(w)}}{1+\norm{w}}
 +\sup_{w\in \ID^{ {p}}}\norm{D^2 g(w)}
+\sup_{w,h\in \ID^{ {p}}}\frac{\norm{D^2 g(w+h)-D^2 f(w)}}{\norm{h}}.
}
Note that for $g\in M'$, $\norm{g}_M \leq \norm{g}_{M'}$.
Defining $I_t(s):=\bone[s\geq t]$, we also typically assume that a test function~$g$ satisfies the smoothness condition that, 
for any $r,s,t$ in $[0,\hT]$  {and $x_1, x_2\in \IR^p$}, 
\ben{\label{eq:smoothgginf}
   \sup_{w\in\ID^{ {p}}}\babs{D^2g(w)\cls{x_1 I_r, x_2(I_s-I_t)}} \Le S_g \, \abs{x_1} \abs{x_2}  \abs{s-t}^{1/2}, 
}
where $S_g$ is some constant depending on $g$  {and $\abs{\,\cdot\,}$ denotes Euclidean norm.}


\subsection{${\rm M/GI}/\infty$ queue}\label{sec:mginfres}
Let $M_n$ be a Poisson process on $\cS:=[0,\hT]\times\IR_+$ with intensity measure
\be{
   \ell_n(dt,dy) \Def n  \alpha(dt)   G(dy),
}
where $\alpha$ is a finite measure on $[0,\hT]$ and $G$ is a distribution function   supported on~$\IR_+$. Viewing~$M_n$ 
as a measure, let $( Y_i,\,i\geq1)$ be i.i.d.\ with distribution function~$\tG$ supported on~$\IR_+$, 
and set 
\[
    N_n \Def M_n + \sum_{i=1}^{x_n} \delta_{(0, Y_i)},
\]
where $x_n\geq0$ is an integer. Then~$N_n$ is a point process with mean measure
\be{
    \lambda_n \Def \ell_n + x_n \bclr{ \delta_0 \times \tG}.
}
For
$(t,y)\in \cS$,  we define 
\be{
    J_{t,y}(s) \Def  \bone\{t\leq s < t+y\}
}
and the  process $X_n:[0,\hT]\to\IR_+$ by
\be{
     X_n(s) \Def \int_{\cS}   J_{t,y}(s) N_n(dt,dy).
}
The process $X_n$ can be regarded as the number of customers in an ${\rm M/GI}/\infty$ queue: A point $(t,y)\in N_n$ 
represents a customer arriving at time $t$ with service time $y$, and such a customer will be in the system at any 
instant  $s$ satisfying $t\leq s < t+y$. We allow the customers initially in the system to have a different service 
distribution, to model the situation where the process is first observed at a typical
time;  {in such a case, the residual service times would have a distribution derived from~$G$, 
but not necessarily the same as~$G$.} 
We consider the ``\emph{heavy traffic}'' regime in which~$n$ is large, so that the total 
rate of arrivals $n\alpha(dt)$ is   large.

Define the centered and scaled random measure 
\be{
   \tN_n \Def \sigma_n^{-1}(N_n-\lambda_n) \quad \text{with $\sigma_n^2 \Eq n$,}
} 
and the process $\tX_n$ by
\ben{\label{eq:mginfqxnt}
     \tX_n(s) \Def \int_{\cS}  J_{t,y}(s) \tN_n(dt,dy).
}
Our main result in this setting is the following.

\begin{theorem}\label{thm:mginfqstn}
 Assume that the convolution $G*\a$ defined by $(G*\a)(s) := \int_0^s G(s-u)\a(du)$, the cumulative 
intensity~$A$ defined by $A(s) := \int_0^s\a(du)$, and the distribution function~$\tG$ are all H\"older continuous 
with exponent~$\b \in (1/2,1]$ and constants $c_{G,\a}$, $c_\a$ and~$c_{\tG}$, respectively.
 Define the measure on $\cS:=[0,\hT]\times\IR_+$ by
\be{
  \Lambda \Def \bclr{ \alpha \times G }+ x \bclr{ \delta_0 \times \tG}, 
}
where $x\geq0$ is fixed. Now set   $\hJ_{(t,y)} := J_{t,y}$ for $t>0$ and $\hJ_{(0,y)} := J_{0,y}-(1-\tG)$,  and let
$Z$ be a   real centred Gaussian process with covariance function, for $0\leq s_1\leq s_2\leq \hT$, 
\ban{
\IE\bcls{Z(s_1) Z(s_2)}
	& =\ \int_{\cS} \hJ_{(t,y)}(s_1) \hJ_{(t,y)}(s_2) \Lambda(dt, dy) \notag \\
	&=\ \int_0^{s_1} \bclr{1-G(s_2-t)} \alpha(dt) 
		+x \, \tG(s_1)\bclr{1- \tG(s_2)}. \label{eq:minfzcov}
}
Then, for any $g\in M'$ either satisfying 
\eq{eq:smoothgginf}, or of the form $g(w)=F\bclr{w(t_1),\ldots, w(t_k)}$ for some {\rm(}twice differentiable{\rm)} 
$F\colon \IR^k \to \IR$ and distinct instants  $t_1,\ldots,t_k\in[0,\hT]$, 
{defining
\ben{\label{ADB-psi-def}
      \ppsi_n(x,x_n,\a,\hT) \Def  
               3\,\frac{\sqrt{\pi x_n/2}  + |x_n -nx|}{2n} + \frac{\alpha\bigl([0,\hT]\bigr)+n^{-1} x_n}{2\sqrt{n}},
}
we have
\ben{\label{eq:mginftfbd}
 \Bigl\vert \mathbb{E}\big[g(\tX_n)\big] -\mathbb{E}[g(Z)] \Bigr\vert \Le 2^{3/2}\|g\|_{M'} { \ppsi_n} (x,x_n,\a,\hT).
}
}
If~$\b = 1$ and $1 \le \a([0,\hT]) \le \a_*\hT$, for some $\a_* < \infty$, and if $\hT \ge n^{1/2}|x_nn^{-1} - x|$,
then, for any $\chi > 0$, there is a constant $K_{\chi}$ such that
\[
    \dlp\bclr{\law(\tX_n), \law(Z)}  \Le  K_{\chi} n^{\chi} \hT^{2/5}  n^{-1/20}. 
\]
where $\dlp$ denotes the L\'evy--Prokhorov distance {\rm(}with respect to Skorokhod topology{\rm)}.
\end{theorem}

\begin{remark}
(1) The constant $S_g$ from the smoothness condition~\eq{eq:smoothgginf} does not appear in the bound, since 
the condition is only used to apply a technical result; see Lemma~\ref{lem:intlin} below.

 \, \,

\noindent (2) As we only consider time intervals $[0,\hT]$, we only require the H\"older continuity on $[0,\hT]$.

\, \,

\noindent (3)  Bounds can also be derived under more general assumptions on~$\a$ and~$\tbe$.  
These can be deduced from the proof of the theorem.
\end{remark}

\begin{remark}
That $\law(\tX_n) \to \law(Z)$ with respect to Skorokhod topology for $x_n=0$ is due to \cite{Borovkov1967} --- 
see the discussion in \cite{Whitt1982} --- and, for general~$x_n$, follows from the results of \cite{Krichagina1997}. 
The only rates of convergence we are aware of are those of  \cite{Besancon2020} and  \cite{Besancon2021},
which give Wasserstein bounds (with respect to the supremum metric) in the special case of the ${\rm M/M}/\infty$ queue. 
Using \cite[Theorem~1.1]{Barbour2021} with the bounds of this paper would lead to rates 
in the Wasserstein distance that would be worse than those derived in \cite{Besancon2021} in this case. 
 {However, our results apply more generally to ${\rm M/GI}/\infty$ queues, which do not appear to be within the scope
of their methods.} 

In the best possible case, where $\beta=1$, our rate of 
convergence for the L\'evy--Prokhorov metric, of $ {\bigo(\hT^{2/5} n^{-1/20+\chi})}$ for any $\chi > 0$,
 is,  to the best of our knowledge, the first rate of convergence in this metric.
A further advantage of our bounds is that the dependence on $\hT$ is explicit, and that $\hT$ could grow like a 
small power of $n$ while still yielding a small bound; this would cover transient approximation almost to stationarity.   
\end{remark}

\begin{remark}\label{rem:mginfgpform}
We can represent $Z$ as a sum of   three independent   centred   Gaussian processes 
$Z= Z_1 + Z_2+Z_3$, where
 \begin{itemize}
 \item[(1)]  $Z_1$ represents the randomness from the services and 
 has covariance structure 
\ba{
\IE\bcls{Z_1(s_1) Z_1(s_2)}
	&\Eq \int_0^{s_1} G(s_1-t)\bclr{1-G(s_2-t)} \alpha(dt),
}
 for
$0\leq s_1\leq s_2\leq \hT$.  
 
\item[(2)]  The second  process $Z_2$ represents the randomness from the arrival process (a kind of weighted 
renewal functional CLT) and is  given by the stochastic integral 
\be{
     Z_2(s) \Eq \int_0^{s} \bclr{1-G(s-t)} B(dt),   
} 
where $B(\cdot)$ is a Gaussian random measure with intensity measure $\alpha(dt)$, that is,    
\ba{
\IE\bcls{Z_2(s_1) Z_2(s_2)}
	&\Eq \int_0^{s_1} \bclr{1-G(s_1-t)}\bclr{1-G(s_2-t)} \alpha(dt); 
}  for
$0\leq s_1\leq s_2\leq \hT$.
\item[(3)] The third process $Z_3$ is a time-changed Brownian bridge: 
  $Z_3(t)= \sqrt{x} B^{br}(\tG(t))$, where $B^{br}$ is a Brownian bridge  with $B^{br}(0) = B^{br}(1)=0$.
\end{itemize} 
For $x=0$, this decomposition is identified in \cite{Borovkov1967}; see also \cite[(2.5) and (2.6)]{Whitt1982} 
and the discussion there. The addition of~$Z_3$ is due to the presence of  customers initially in the system. 
 {The number of those remaining in the system at time~$t$ is just the number with service time greater than~$t$, 
the empirical complementary cumulative distribution function at~$t$.  This,} after scaling,  
converges as a function of~$t$ to a time--changed Brownian bridge. 
\end{remark}

\begin{remark}\label{rem:fddsmmi}
If $g(w)=F\bclr{w(t_1),\ldots, w(t_k)}$, and $F$ is bounded with bounded partial derivatives of order up to three, 
then $\norm{g}_M\leq  c k^3$, where~$c$ is an upper bound for~$F$ and its first three partial derivatives. Assuming 
that $\abs{x_n- nx} =\bigo(n^{1/2})$, the bound~\eqref{eq:mginftfbd} is of order $\bigo(k^3 n^{-1/2})$, and standard smoothing 
arguments, as in \cite[Section~1.1.4]{Gotze2019}, imply that for any convex $\cK\subseteq\IR^k$ and $\eps>0$, 
\be{
   \bbabs{\IP\bclr{(\widetilde X_n(t_i))_{i=1}^k\in \cK} - \IP\bclr{(Z(t_i))_{i=1}^k}\in \cK} 
           \leq C\bclc{ \eps^{-3} k^3  n^{-1/2} + k^{1/4} \eps}.
}
Choosing $\eps =k^{11/16} n^{-1/8}$, leads to a uniform upper bound of order $k^{15/16}n^{-1/8}$ on the difference of 
convex set probabilities for \emph{any} $k$-dimensional distributions of $\widetilde X_n$ and $Z$. The power of~$k$ 
in the bound is likely not optimal, but still leads to meaningful results for~$k$ growing like a small power of~$n$.
\end{remark}

%

\subsection{$\text{GI}/\text{GI}/\infty$ queue}\label{sec:mggres}

Consider a stationary renewal process $V_n$ on $[0,\hT]$,  whose renewal distribution~$\nu_n$ is that of~$R/n$, where~$R$ is
a positive integer valued random variable with  \emph{aperiodic} 
support having mean~$m$, variance~$v^2$ and $\IE[R^r]<\infty$ for some fixed $r\geq 5$.
Let~$G$ be a distribution function supported on $\IR_+$ and let $(Y_i,\,i\geq1)$ be an i.i.d. sequence with 
distribution~$G$ that is also independent of the  renewal process~$V_n$. 
Now define  the random measure~$N_n$ on $\cS:=[0,\hT]\times \IR_+$ by
\be{
   N_n \Def \sum_{i=1}^{\floor{n\hT}} M_n(i/n) \delta_{(i/n, Y_i)},
}
where now
\[
    M_n(s) \Def  V_n(s) - V_n(s-) \Def  \mathbf{1}[\text{a renewal occurs at time $s]$, \quad   $s\in[0,\hT]$},  
\] 
and~$Y_i$ represents the service time of a customer arriving at time $i/n$. 
Due to the stationarity, the mean measure of~$N_n$ is 
\be{
    \lambda_n \Def \frac{1}{m} \sum_{i=1}^{\floor{n\hT}} \bclr{ \delta_{i/n} \times G}. 
}
For $(t,y)\in \cS$, we define $ J_{t,y}(s) := \mathbf{1}[t\leq s < t+y]$
and  
\be{
    X_n(s) \Def \int_{\cS}   J_{t,y}(s) N_n(dt,dy).
}
As for the ${\rm M/GI/\infty}$ queue, the process $X_n$ can be regarded as the number of customers in a 
$\text{GI}/\text{GI}/\infty$ queue with arrival times given by a stationary renewal process driven by~$\nu_n$ 
and service times distributed according to~$G$. 
Note that because $\nu_n$ (the law of $R/n$) is discrete, the results of the previous section are not a special 
case of those derived here.

Define the centered and scaled random measure 
\be{
    \tN_n \Def \sigma_n^{-1}(N_n-\lambda_n),
}
where 
\ben{\label{ADB-variance-def}
  \sigma_n^2 \Def nv^2/m^3,
}
and define the process $\tX_n$ by
\be{
  \tX_n(s) \Def \int_{\cS}  J_{t,y}(s) \tN_n(dt,dy).
}
Writing $I_t(\cdot):=\bone\{ \cdot \geq t\}$ and  $\overline G:= 1-G$,
our main result of the section is the following. 

\begin{theorem}\label{thm:gginfqstn}
Recall the notation just above.  Assume that $\ex R^5 < \infty$ and that, 
for some   $\beta\in(0,1]$ and $0 < \LL < 1$, the distribution function $G$ satisfies
\ben{\label{eq:gcon}
  G(t) - G(s) \Le g_G(s)(t-s)^\b, \quad 0 < s < t \le \hT,   
}
where $g_G\colon \IR_+ \to \IR_+$ is bounded and non-increasing and such that $\int_0^\infty g_G^\LL(s)\,ds < \infty$.
Define the measure on $\cS:=[0,\hT]\times\IR_+$ by
\be{
   \Lambda(dt,dy)\Def dt\, G(dy),
}
set 
\[
    \hJ_{t,y}(s) \Def \frac{m}{v}\,J_{t,y}(s)-\frac{m+v}{v}\,\overline G(s-t) I_t(s),
\] 
and let~$Z$ be a Gaussian process with covariance function, for $0\leq s_1\leq s_2\leq \hT$,
\ban{
\IE\bcls{Z(s_1) Z(s_2)}\ 
&=\ \int_{\cS} \hJ_{t,y}(s_1) \hJ_{t,y}(s_2) \Lambda(dt, dy) \notag \\
	&=\ \frac{m^2}{v^2} \int_{0}^{s_1} \overline G(s_2-t)  G(s_1-t) dt +  \int_{0}^{s_1} \overline G(s_1-t) 
                      \overline G(s_2-t) dt.\label{eq:covgginfz}
}
Then, for any $g\in M'$ satisfying the smoothness condition~\eq{eq:smoothgginf} for some $S_g$, there is a 
constant $C$ depending on $\law(R)$  and $g_G(0)$ such that,  with $\bab := \min\{\beta,1/2\}$, we have  
\ban{
\abs{\IE[g(&\tX_n)] -\IE[g(Z)]}  \Le C\hT\bclr{S_g n^{-1/2} +  \norm{g}_{M'} n^{-\bab}}. 
            \label{GGI-expec-bnd}
}
If $g\in M'$ is  of the form $g(w)=F\bclr{w(t_1),\ldots, w(t_k)}$ for some {\rm(}twice differentiable{\rm)} 
$F\colon \IR^k \to \IR$ and  distinct instants  $t_1,\ldots,t_k\in[0,1]$, 
then the same bound holds, but with $S_g n^{-1/2}$ replaced by $\norm{g}_{M'} k^3 n^{-1}$.

Moreover, assuming that $\hT \ge 1$,
\ignore{
\[
     x(\z) \Def \min\{1/2,\z(r-2),\b - \z(1+\b)\},
\]
}
taking $l_r := \lceil (r - \LL(r-1))/2 \rceil - 1$, $\b_r := \b(r-2)/(r-1)$, we have
\ben{\label{eq:gginfrate}
    \dlp\bclr{\law(\tX_n), \law(Z)} 
              \Le C \sqrt{\log n}  \Bigl\{(\hT^{4}n^{-\bab})^{(l_r\b_r-1)}T^{3}\Bigr\}^{1/(6l_r + 4l_r\b_r -1)} ,
} 
\ignore{
where~$\z$ is the unique solution in $(0,1/2]$ to the equation
\[
        x(\z) \Eq \frac{ - (6l_r + 4l_r\b_r -1)(2\z - 1/2) + \tfrac18(32l_r + 28l_r\b_r +3)\log_n\hT }
                       { 4l_r + 3l_r\b_r - 1 }\,,
\]
and 
}
where $\dlp$ denotes the L\'evy--Prokhorov distance  {\rm(}with respect to the Skorokhod topology{\rm)}.
For instance, if the distribution of~$R$ has all positive moments, and if~$G$ has a finite moment and a bounded 
and ultimately monotone density, and if~$\hT \le n^{\psi}$ for some $\psi < 1/8$, then, for any $\chi > 0$,
there is a constant~$K_\chi$ such that
\[
      \dlp\bclr{\law(\tX_n), \law(Z)}  \Le K_\chi n^\chi  \hT^{2/5} n^{-1/20}.
\]

\end{theorem}

\ignore{
Moreover, assuming $\hT\geq 1$, if in~\eq{eq:gcon} we can take $\beta\in(1/2,1]$, then with $\dlp$ denoting the 
L\'evy--Prokhorov distance  {\rm(}with respect to Skorokhod topology{\rm)}, there is are constants $C$ and $\iota$ 
depending on $\law(R), r$ and $c$ such that for any $\mu\in \bclr{0,\frac{(2\beta-1)(2r-1)}{2(\beta(r-1)-1)}}$, 
\ben{\label{eq:gginfrate}
\dlp\bclr{\law(\tX_n), \law(Z)} \leq C \hT^{(2r+3)/2} 
             \bbclr{\sqrt{\log(n)} n^{-\phi(\mu)}+ ne^{- \iota \hT n^{1/2-\mu (2\beta-1)/5}}},
}
where $\phi(\mu)$ can be taken to be equal to 
\be{
\min\bclc{ \nu, (r-2)\ee-3(\mu+\nu), \beta - \ee(1+\beta) - 3(\mu+\nu), \tsfrac{1}{2}- 3(\mu + \nu), 
                  \tsfrac{1}{2} - 2(\mu +\nu + \ee), \mu \tsfrac{\beta(r-1)-1}{2r-1}},
}
for any $\nu>0$ and $0< \ee < \beta/(1+\beta)$. In particular, setting $\mu=\frac{5}{18(2\beta + 4)}$, $\phi$ 
can be taken to be at least as large as $\frac{2\beta -1}{18 (2\beta +4)}$ for any $r\geq 5$.
}

\begin{remark}
The expression~\eq{eq:covgginfz} for the covariance agrees (up to scaling) with \cite[(2.5) and~(2.7)]{Whitt1982}. 
Note also that choosing~$G$ supported on $(\hT,\infty)$ corresponds to the renewal CLT, where the limiting Gaussian process 
is the standard  Brownian motion. The decomposition in Remark~\ref{rem:mginfgpform} still applies here with 
$Z=(m^2/v^2) Z_1 + Z_2$ and $\alpha(dt)=dt$. The addition of initial customers would contribute to the limit in the 
same way as in Theorem~\ref{thm:mginfqstn}; that is, it would add a Brownian bridge component as described in 
Remark~\ref{rem:mginfgpform}, but we omit this for the sake of clarity. Weak convergence with respect to the Skorokhod 
topology was (essentially) shown in \cite{Borovkov1967}, and here we can view~\eq{eq:gginfrate} as a rate of convergence. 
We are not aware of previous results with rates of convergence.  The same considerations as  in Remark~\ref{rem:fddsmmi} 
lead to a uniform upper bound of order $k^{15/16}n^{-\bab/4}$ on the difference of convex set probabilities 
for \emph{any} $k$-dimensional distributions of $\widetilde X_n$ and $Z$.   Again, we remark that an appealing 
aspect of our 
bound is the explicit incorporation of~$\hT$.
\end{remark}

\begin{remark}
For a process $ {X_n'}$ defined analogously to~$X_n$, but driven by~$V_n'$ defined to be a delayed 
(rather than stationary) 
renewal process with inter-renewal distribution $\law(R/n)$, it is easy to see that 
 {$X_n$ and $X_n'$ can be constructed on the same space so that} $\norm{X_n- {X_n'}}$ is 
stochastically dominated by the coupling time 
\be{
   T_c \Def \inf\bclc{i\geq 1\colon M_n(i/n) = M_n'(i/n)=1},
} 
where $M_n'(s) = V_n'(s) - V_n'(s-)$ is defined in analogy with $M_n$ (here, we view $V_n$ and~$V_n'$ as 
renewal processes 
on $[0,\infty)$, and so $T_c$ does not depend on~$n$). Defining the scaled process 
\be{
  {  \widetilde X_n'} \Def \frac{ {X_n'} - \int_{\cS} J_{t,y} \lambda_n(dt,dy)}{\sigma_n},
}
we thus easily find that $\norm{\tX_n -  {\tX_n'}}$ is stochastically dominated by  $\bclr{m^{3/2}T_c/v} n^{-1/2}$. 
Under the hypotheses of the theorem, \cite[Proposition (6.10)]{Pitman1974} implies that~$T_c$ is finite 
with probability one, 
and if the delay distribution has finite $(r+1)$-moment, that $\IE[T_c^{{r}}]<\infty$. 
Under this moment condition, 
for any function $g\in M'$ with $\norm{Dg(w)}<\infty$, we then easily have that
\be{
    \babs{\IE[g( {\tX_n'})] - \IE[g(\tX_n)]} \Le \norm{Dg(w)} \frac{m^{3/2} \IE[T_c]}{v \sqrt{n}}.
}
Combining this with the bounds of the theorem and the triangle inequality gives bounds on 
$\babs{\IE[g( {\tX'_n})] - \IE[g(Z)]}$ for $g\in M'$ with  $\norm{Dg(w)}<\infty$, and subsequent bounds of the same 
order on $\dlp$ (the restriction that $\norm{Dg(w)}<\infty$ is no problem, since the key result 
Theorem~\ref{thm:smooth} only requires bounds on a smaller class of test functions with bounded derivatives).
\end{remark}

\subsection{General approximation theorem}

Here we state the general approximation theorem used to prove Theorems~\ref{thm:mginfqstn} and~\ref{thm:gginfqstn}.  
We first need a definition.

\begin{definition}\label{def:redpalm}
For a point process $\Xi\subseteq \IR^d$ with mean measure $\kappa$, we say $\Xi\sp{u}$ is distributed as the 
reduced Palm measure of $\Xi$, if 
\be{
  \IE \bbcls{ \int g\bclr{\Xi, u}\Xi(du) } \Eq \IE \bbcls{\int g\bclr{\Xi\sp{u}+\delta_u, u} \kappa(du)}
}
for all functions $g$ such that the integral on the left hand side exists.
\end{definition}
For simple point processes (meaning there is at most one point at any location),  we think of $\Xi\sp{u}+\delta_u$ 
as having the distribution of $\Xi$ conditional on there being a point at $u$, which explains why we can 
take $\Xi\sp{u}=\Xi$ if $\Xi$ is a Poisson process. 
For rigorous background on reduced Palm measures, see \cite[Chapter~13]{Daley2008}. 

\begin{theorem}\label{thm:stnpalm}
Recall the notation and definitions of Section~\ref{sec:setup} leading up to~\eq{eq:gentildex}. 
Let $(N, N\sp{u})_{u\in\IR^d}$ be a collection of couplings of $N$ with its reduced Palm measures and define 
\[
    X\sp{u}(s) \Def \int_{\IR^d}J_{v}(s) N\sp{u}(dv), ~s\in[0,\hT].
\]
Let $Z$ be a centered Gaussian process with  almost surely continuous sample paths with
\ben{\label{eq:gpcovexp}
  \IE[Z(s)Z(t) {^\tr}] \Eq \int_{u\in \IR^d} \hJ_u(s)\hJ_u(t) {^\tr} \Lambda(du),
} 
where the function $(u,s)\longmapsto \hJ_u(s)\in\mathbb{R} {^p}$ is measurable such that 
\begin{center}
    $\hJ_u(\cdot)\in L^2([0,\hT] {\to\IR^p})$  and  $\hJ_{\cdot}(s)\in L^2(\IR^d {\to \IR^p}, \Lambda)$, 
          for all $s\in[0,\hT]$ and all $u\in\mathbb{R}^d$. 
\end{center}
Now suppose that $g\in M$, and define  $f:=f_g$ to be the Stein solution given in Theorem~\ref{prop:stnsol}. 
Then  we have
\ban{
\abs{&\IE[g(\tX)] -\IE[g(Z)]} \notag\\
	&\leq\ \bbbabs{\IE \bbbcls{ \int_{\IR^d} D^2 f(\tX)[\hJ_{u}\uts] \Lambda(du)
    - \int_{\IR^d} D^2 f(\tX)\bcls{J_{u},\,\IE[X\sp{u}-X+J_u]}\bclr{\sigma^{-2}\lambda(du)}}} \label{eq:stnbdt1} \\
    &\qquad + \bbbabs{\IE \bbbcls{\int_{\IR^d} 
            D^2 f(\tX)\bcls{J_u,\, (X\sp{u}-X)-\IE[X\sp{u}-X]} \bclr{\sigma^{-2}\lambda(du)}}}\label{eq:stnbdt2} \\
   &\qquad + \frac{\norm{g}_M}{2\sigma}
    \IE \bbbcls{\int_{\IR^d}\norm{J_{u}}\norm{X\sp{u}-X+J_u}^2\bclr{\sigma^{-2}\lambda(du)}}. \label{eq:stnbdt3}
}
If $N$ is a Poisson point process, $\Lambda = \sigma^{-2} \lambda$ and $\hJ_u = J_u$, then
for any $g\in M$, 
\ben{\label{eq:poisstnbd}
   \abs{\IE[g(\tX)] -\IE[g(Z)]} 
	\Le  \frac{\norm{g}_M}{2\sigma} \bbbcls{\int_{\IR^d}\norm{J_{u}}^3\bclr{\sigma^{-2}\lambda(du)}}.
}
\end{theorem}

\begin{remark}\label{rem:checkz}
(1) Although our focus in this work is about real Gaussian process approximation,  there is no essential extra difficulty in establishing the above multivariate Gaussian approximation result.  See   Section \ref{sec:mginf}  for an illustration in the multivariate setting. 

 \, \,

\noindent 
(2)
To check the hypothesis that~$Z$ has continuous sample paths,  it is enough to establish that, for some positive 
constants $C,b$, we have
\ben{\label{eq:gpkolcont}
   \IE\bcls{\babs{Z(s)-Z(t)}^2} \Le C \abs{s-t}^b,
}
which, using Gaussianity, implies the Kolmogorov continuity criterion.

 \, \,

\noindent(3) In general, the terms in  \eqref{eq:stnbdt3}
 and  \eqref{eq:poisstnbd} are not automatically finite. In our applications, the integrand $J_u$ is uniformly bounded, $\IE[ \|X\|^2] $ is finite and $N$ is a simple point process over a subset $\mathcal{S}$ of $\IR^d$ with finite  intensity measure such that $N(\mathcal{S})$ has finite second moment. Under these extra assumptions, the aforementioned terms are  finite. 
 
\end{remark}
\begin{remark}
In practice, the choice of $\Lambda$ and $\{\hJ_u\}_{u\in\IR^d}$ arises from computing the quantity
$\IE[X\sp{u}-X+J_u]$, plugging the resulting expression into the second  term of the difference of~\eq{eq:stnbdt1}, 
and then discarding the asymptotically negligible terms. This is an appealing aspect of the theorem, as it suggests 
a candidate limit, while also providing an intuitive expression for the covariance of the limit in the 
form~\eq{eq:gpcovexp}.
\end{remark}

\begin{remark}
Bounding~\eq{eq:stnbdt1} and~\eq{eq:stnbdt2}  in applications requires using the structure of $(X\sp{u}-X)$ 
and its mean, 
along with the bounds and ``smoothness'' properties of $f$ given in Theorem~\ref{prop:stnsol} below. 
Bounding~\eq{eq:stnbdt2} is the main difficulty in applying the theorem, and typically requires constructing 
intermediate couplings that exploit local or weak global dependence.
\end{remark}

To apply the theorem for the ${\rm M/GI/\infty}$ queue, we need to define the reduced Palm couplings of the 
arrival/service point process.  Away from time zero, the arrivals/services are a Poisson point process, and 
so we can take the reduced Palm measure to be the original process. For the customers in the system at time zero, 
the Palm measure corresponds to removing a point at random, which is only a small 
perturbation of the original process. 
For the $\text{GI}/\text{GI}/\infty$ queue, the arrival/service point process is no longer Poisson. 
However, because the service 
times are i.i.d., constructing the reduced Palm measure coupling at a point $(s,y)$  comes down to constructing a 
close coupling of a stationary renewal process to one conditioned to have a renewal at $s$, which in turn is similar 
to coupling a stationary renewal sequence to a zero-delayed renewal process, and this is well-understood. The details 
are in Sections~\ref{sec:mginf} and~\ref{sec:gginf}.

Theorem \ref{thm:stnpalm} follows from a new development of Stein's method \cite{Stein1972, Stein1986},
formulated as Theorem~\ref{prop:stnsol}. 
Stein's method provides a general framework for bounding the error when approximating a complicated distribution 
of interest by a well-understood target distribution; see \cite{Ross2011} for a basic introduction. By now, Stein's 
method has been developed  for a large number of univariate distributions, as in the monographs  \cite{Chen2011} 
for the normal and \cite{Barbour1992} for the Poisson. Stein's method for multivariate distributions other than the 
normal is not so well developed, and even less is known for random processes.  
 {Poisson process approximation is a notable exception, with a succession of papers going back to  
\cite{Barbour1988} and \cite{Arratia1989}.  There is also work on approximation by Brownian motion, which began 
with \cite{Barbour1990}, }
and on some closely related Gaussian processes, such as time changes of Brownian motion \cite{Kasprzak2017a, Kasprzak2020a} 
and multivariate correlated Brownian  motions \cite{Kasprzak2020, Dobler2021}; there has also been recent work for 
Dirichlet Process approximation \cite{Gan2021}.

All of the Gaussian process approximation results just cited are derived using Barbour's generator 
approach \cite{Barbour1990}, which identifies a ``characterizing'' operator of a Gaussian process as the 
generator of the Ornstein--Uhlenbeck semigroup. To avoid working with the generator in the continuum, these 
papers first approximate the process of interest by a discretized version of the Gaussian process, which has 
an Ornstein--Uhlenbeck generator of a simple form. After this is achieved, the problem is reduced to showing that
the discretized Gaussian process is close to the true Gaussian process.  Here we avoid this two step procedure, 
by developing the relevant properties of the ``Stein solution'' for any Gaussian process with  continuous 
sample paths; see Theorem~\ref{prop:stnsol} below. This theorem can be used to 
prove approximation bounds 
in quite general settings, in which the dependency structures are amenable to Stein's method, 
such as those exhibiting an exchangeable pair. The formulation is particularly
useful in our applications,
where the jumps of the processes that we study occur at 
random times. A final remark on this point is that any disretization error between the process and the target 
Gaussian process is captured by~\eq{eq:stnbdt1} from Theorem~\ref{thm:stnpalm}, and bounding this term typically 
relies on the smoothness property~\eq{eq:smoothgginf}.

The recent papers \cite{Coutin2020a} and \cite{Besancon2021} use Stein's method 
to obtain bounds  {in the bounded Wasserstein distance for Donsker's theorem, and for Lipschitz functionals of Poisson 
measures. 
As discussed in \cite{Barbour2021}, they obtain rates of convergence in this restricted setting of better order than those that
our method typically yields, 
though less good than those obtainable using strong approximation.  However, in their approach, they make use of the independence
structure within the process being approximated, and of the fact that the limiting process is Brownian motion; 
we need neither of these simplifications.}

There is also an approach to Stein's method on Hilbert and abstract Wiener spaces, initiated in
\cite{Shih2011} and developed further in \cite{Coutin2013, Besancon2020, Bourguin2020} and \cite{Bourguin2021}.   
These approximations are not sufficient to directly imply weak convergence with respect to the Skorokhod topology,
and hence do not imply rates for such convergence, either; 
see the nice discussion of some of these approaches in \cite[Section~1.5]{Dobler2021}.

The remainder of the paper is organised as follows. In Section~\ref{sec:stnmetgp}, we develop Stein's method
in the general context of Gaussian process approximation, establishing Theorem~\ref{prop:stnsol},
together with some ancillary results.  We then
prove Theorem~\ref{thm:stnpalm}. In Sections~\ref{sec:mginf} and~\ref{sec:gginf} 
we apply Theorem~\ref{thm:stnpalm} to the ${\rm M/GI/\infty}$ and $\text{GI}/\text{GI}/\infty$ queue examples given above, 
proving Theorems~\ref{thm:mginfqstn} and~\ref{thm:gginfqstn}.

\section{Stein's method for Gaussian processes}\label{sec:stnmetgp}

Our first step in developing Stein's method for Gaussian processes is to establish a useful 
form for the characterizing operator.

\begin{proposition}\label{prop:stncharop}
Let $Z$ be a real centred continuous  $\IR^p$-valued Gaussian process on $[0,\hT]$ with 
\ben{\label{eq:covform}
    K(s_1,s_2) \Def \IE[Z(s_1)Z(s_2) {^\tr}] \Eq \int_{\IR^d} J_u(s_1)J_u(s_2) {^\tr} \Lambda(du),
}
where~$\Lambda$ is a measure on $\IR^d$ 
and the function $(u,s)\longmapsto   J_u(s)\in\mathbb{R} {^p}$ is measurable such that 
\begin{center}
   $J_u(\cdot)\in L^2([0,\hT] \to \IR^p)$    and   $J_{\cdot}(s)\in L^2(\IR^d\to \IR^p, \Lambda)$, 
        for all $s\in[0,\hT]$ and  all $u\in\mathbb{R}^d$.
\end{center}
Then
for any function $f\in M$,
\ben{\label{eq:2ordt}
\IE \Big(  D^2 f(w)[Z\uts] \Big) \Eq \int_{\IR^d} D^2 f(w)[ J_u\uts]\Lambda(du),
}
and
\ben{\label{eq:charop1}
 \IE\left( \int_{\IR^d} D^2 f(Z)[ J_u\uts]\Lambda(du)- D f(Z)[Z] \right) \Eq 0.
}
\end{proposition}

\begin{proof}
First we show that all the expectations exist. Since $f\in M$, we have
\ba{
\abs{D^2f(w)[Z\uts]} \Le  \norm{D^2f} (1+\norm{w})\norm{Z}^2 \quad{\rm and}\quad 
\abs{D f(Z)[Z]}  \Le  \norm{D f}(1+\norm{Z}^2)  \norm{Z},
}
so that we need to show that $\norm{Z}$ has finite third moment,  which is guaranteed by Fernique's theorem.
To establish the expressions for the moments, we use the  {multivariate}  Karhunen--Lo\`eve expansion of~$Z$; 
see  {\cite[Section~2.2]{Happ2018}}. 
Define the linear operator $T$ {on the Hilbert space   $L^2([0,\hT] {\to \IR^p}) $} by setting
\[
(Tf)(s)  \Def  \int_0^\hT  K(s,t) f(t) dt.
\]
{It is easy to see that $T$ is a  positive and compact self-adjoint operator on  $L^2([0,\hT] {\to \IR^p})$, 
so that by the spectral theorem, we can find  $\{h_k, k\in\N\}$ that is an orthonormal basis of $L^2([0, \hT] {\to \IR^p})$ 
formed by the eigenvectors  of $T$ with respective eigenvalues $\{\ell_k, k\in\N\}\subset\R_+$. 
As a result,}
 $Z(t)$ admits the following representation
\begin{align}
Z(t)  \Eq  \sum_{k\in\N} X_k h_k(t),  \label{USE}
\end{align}
where
\[
X_k  \Def  \int_0^\hT Z(t) {^\tr}  h_k(t)\, dt \,  { = \int_0^\hT   h_k(t)^\tr Z(t) \, dt \,,}
\]
and the convergence in~\eq{USE} can be taken in $L^2(\Omega)$, uniformly in $t\in[0,\hT]$, and, because of the 
assumption of continuity of sample paths, can also be taken with respect to sup norm; see \cite[Theorem~3.1.2]{Adler2007}.

For  $k,j\in\N$, we have 
\begin{align*}
  \E\big[ X_k X_j \big] &\Eq \E\left[  \int_0^\hT  \int_0^\hT { h_j(s)^\tr Z(s)Z(t)^\tr \,  \, h_k(t)\,} dt ds \right] 
                 \Eq \int_0^\hT  \int_0^\hT { h_j(s)^\tr K(s,t) \, h_k(t)}\, dt ds  \\
    &\Eq \int_0^\hT  {h_j(s)^\tr} (Th_k)(s)  \, ds  \Eq   \ell_k \mathbf{1}[k=j],
\end{align*}
and so the variables $\{X_k\}_{k\in \IN}$ are independent centred Gaussian random variables.

Using the representation~\eqref{USE}, we have,  for the symmetric bilinear form~$A = D^2f(w)$, 
\[
  A[Z,Z] \Eq A\left[   \sum_{k\in\N} X_k h_k,\,  \sum_{j\in\N} X_j h_j\right]  \Eq \sum_{k,j\in\N} X_kX_j A[h_k, h_j],
\]
so that 
\ben{\label{eq:expazz}
\E\big[ A[Z,Z] \big]  \Eq \E\left[ \sum_{k,j\in\N} X_kX_j A[h_k, h_j] \right] \Eq \sum_{k\in\N} \ell_k A[h_k, h_k] \,.
}
On the other hand, expanding $J_u$ in the orthonormal basis $\{ h_k\}$ implies 
\[
    J_u  \Eq \sum_{k\in\N}  h_k  J_u\sp{k},
\]
where $J_u\sp{k} :=  \int_0^\hT J_u(s)^\tr h_k(s)\, ds \in\IR$ for each $k\in\N$.
It follows that
\besn{\label{eq:expajj}
 \int_{\R^d} A\big[J_u,  J_u \big] \, \Lambda(du) 
      &\Eq \int_{\R^d} A\left[ \sum_{k\in\N}  h_k  J_u\sp{k} ,\, \sum_{j\in\N}  h_j J_u\sp{j}\right] \, \Lambda(du) \\
      &\Eq \sum_{k,j\in \N}  A[h_k, h_j] \int_{\R^d} J_u\sp{k} J_u\sp{j} \Lambda(du).
}
Continuing with straightforward calculations, we have
\begin{align*}
\int_{\R^d} J_u\sp{k}J_u\sp{j} \Lambda(du)
	& \Eq  \int_{\R^d} \left( \int_0^\hT  h_k(s){^\tr} J_u(s)\, ds \right) \left( \int_0^\hT  J_u(t){^\tr}  h_j(t)\, dt \right)\Lambda(du ) \\
	& \Eq \int_0^\hT\int_0^\hT h_k(s)^\tr \left\{  \int_{\R^d} J_u(s)J_u(t) {^\tr}\Lambda(du ) \right\}  h_j(t) dsdt \\
& \Eq  \int_0^\hT\int_0^\hT   { h_k(s)^\tr}K(s,t) h_j(t) dsdt  \Eq  \int_0^\hT   {h_k(s)^\tr (Th_j)(s) ds}\\
& \Eq  \ell_j \int_0^\hT { h_k(s)^\tr h_j(s)ds}  \Eq  \ell_j \mathbf{1}[ k=j].
\end{align*}
Plugging this into~\eq{eq:expajj}, and noting \eq{eq:expazz} gives~\eq{eq:2ordt}:
\begin{align}
 \int_{\R^d} A\big[J_u,  J_u \big] \, \Lambda(du)  & \Eq  \sum_{k \in\N} \ell_k A[h_k, h_k]  \Eq   \E\big[ A[Z,Z] \big] \,. 
             \label{eq:intmedaj}
\end{align}
For~\eq{eq:charop1}, the first equality in~\eq{eq:intmedaj} implies that it is enough to establish that
\ben{\label{eq:serisow}
\IE\bcls{D f(Z)[Z]} \Eq \sum_{k\in \IN} \ell_k \IE\bcls{D^2 f(Z)[h_k,h_k]}.
}
Writing $Z_k:= Z - X_k h_k$, we have
\be{
\IE\bcls{D f(Z)[Z]}
	 \Eq  \sum_{k\in \IN} \IE \bcls{ X_k D f(Z_k + X_k h_k)[h_k]}
	 \Eq \sum_{k\in \IN} \IE \Big[  \IE \Big( X_k D f(Z_k + X_k h_k)[h_k]  \big\vert Z_k \Big) \Big] .
}
Since $Z_k$ is independent of $X_k$, we can apply the one dimensional Stein identity
\be{
\IE[ X_k g(X_k) |Z_k]  \Eq  \ell_k \IE[g'(X_k) |Z_k]
}
with $g(x)=D f(Z_k + x h_k)[h_k] $ to each term in this sum. Then~\eq{eq:serisow} easily follows by noting 
that $g'(x)= D^2 f(Z_k + x h_k)[h_k, h_k]$ and thus $g'(X_k)= D^2 f(Z)[h_k, h_k]$.
\end{proof}

The next result represents the foundation of Stein's method for continuous Gaussian processes.

\begin{theorem}\label{prop:stnsol}
Let $Z$ be a centred continuous Gaussian process with   covariance function given by \eqref{eq:covform}. 
Given $g\in M$, we define $f_g:\ID {^p}\to\IR$ by 
\ben{\label{eq:stnsol}
f_g(w) \Def -\int_{0}^\infty\bbclr{ \IE\bcls{g(w e^{-s}+ \sqrt{1-e^{-2s}} Z)} - \IE\bcls{ g(Z)}} ds.
}
Then $f_g\in M$  and for $k\in\{1,2\}$,
\ben{\label{eq4.3}
D^k f_g (w)  \Eq  -\IE \int_{0}^\infty e^{-ks}D^k g(w e^{-s}+ \sqrt{1-e^{-2s}} Z) ds.
}
Furthermore, for $w,w',w_1,w_2\in \mathbb{D}^p$, we have
\ban{
\babs{D^2 f_g(w+w')[w_1,w_2]-D^2f_g(w)[w_1,w_2]}
	& \Le    \norm{g}_M  \norm{w_1}\norm{w_2} \norm{w'}, \label{eq:3dstnbd}
}
and if $g\in M'$, then 
\ban{
\babs{D^2 f_g(w)[w_1,w_2]} 
	& \Le  (3/2) \norm{g}_{M'} \norm{w_1}\norm{w_2},\label{eq:2dstnbd} 
}
Finally, $f_g$ satisfies the Stein's equation
\ben{\label{SteinEq}
\cA f_g(w) \Def \int_{\IR^d} D^2 f_g(w)[ J_u\uts]\Lambda(du)-D f_g(w)[w] \Eq g(w) -\IE\bcls{g(Z)}.
}
\end{theorem}

\begin{proof}
That $f_g\in M$ and that~\eq{eq4.3} holds follow from the arguments of 
\cite[Lemma~4.1]{Kasprzak2017} (see also \cite{Barbour1990}) for the special case of Brownian motion. Their 
argument only relies on the supremum of the Gaussian process having finite third moment, which is also valid in our setting. 

The bounds on the derivatives also more or less follow along the same lines as existing work, 
see \cite[Proposition~3.2]{Kasprzak2020a} or \cite[Proposition~5.5]{Kasprzak2020}, but our setting is different 
enough that we include a proof.
To show~\eq{eq:2dstnbd}, we use equation~\eq{eq4.3} and Lemma~\ref{lem:bilnbds} below, which relates the 
absolute value of a bilinear form at a given argument to its norm, to find that
\begin{align*}
\babs{D^2 f_g(w)[w_1,w_2]}
	& \Le  \int_{0}^\infty e^{-2s} \IE\bbcls{\babs{D^2 g(w e^{-s}+ \sqrt{1-e^{-2s}} Z) [w_1,w_2]}} ds \\
	& \Le  3 \norm{w_1}\norm{w_2} \int_{0}^\infty e^{-2s} \IE\bbcls{\bnorm{D^2 g(w e^{-s}+ \sqrt{1-e^{-2s}} Z)} } ds \\
	& \Le  3 \norm{g}_{M'} \norm{w_1}\norm{w_2} \int_{0}^\infty e^{-2s} ds  \Eq (3/2) \norm{g}_{M'} \norm{w_1}\norm{w_2},
\end{align*}
where the third inequality uses that $g\in M'$. The proof of~\eq{eq:3dstnbd}   follows from similar arguments 
using equation \eqref{eq4.3}, Lemma \ref{lem:bilnbds} and the Lipschitz continuity of $D^2g$.

The usual path to show~\eq{SteinEq} is to view the family of operators 
$P_s\colon g\mapsto g(w e^{-s}+ \sqrt{1-e^{-2s}} Z)$ 
as an Ornstein--Uhlenbeck semigroup with generator equal to the characterizing operator~$\cA$, and then the 
result follows essentially from strong continuity of the semigroup. However, the semigroup is not strongly 
continuous, even for~$Z$ a Brownian motion, and so
an alternative approach is to follow the proof of the relevant result for strongly continuous semigroups; 
see  \cite{Kasprzak2017}. While such a  strategy could work in  our setting, we provide a direct proof that is 
simpler than existing approaches, using Gaussian calculations and~\eq{eq4.3}.

Putting $\mathcal{W}_s= w e^{-s} + \sqrt{1-e^{-2s}} Z$ and using~\eq{eq:2ordt}, \eq{eq4.3}, we can rewrite  
the left hand side of~\eqref{SteinEq}  as 
 \begin{align}
 \E\int_0^\infty e^{-s} Dg(  \mathcal{W}_s  )[\w] ds - \E \int_0^\infty e^{-2s} D^2g(  \mathcal{W}_s)[Z', Z'] ds, 
              \label{eq:lhsstneq}
 \end{align}
 where $Z'$ is an independent copy of $Z$.
The right hand side of \eqref{SteinEq} can be written as
 \[
 \E\big[ g(\w) - g(Z) \big]  \Eq  - \E\int_0^\infty \frac{d}{ds} g\big( \mathcal{W}_s \big) ds,
 \]
and, since $g\in M$,
\[
\frac{d}{ds} g\big( \mathcal{W}_s \big)  \Eq   Dg(\mathcal{W}_s)\Big[ -e^{-s} \w + \frac{e^{-2s}}{\sqrt{1- e^{-2s}}} Z \Big].
\]
Thus, using linearity of the derivative,
\begin{align*}
g(\w) - \IE\bcls{g(Z)}
   &\Eq  - \E\int_0^\infty    Dg(\mathcal{W}_s)\Big[ -e^{-s} \w + \frac{e^{-2s}}{\sqrt{1- e^{-2s}}} Z \Big] ds \\
  &\Eq \E\int_0^\infty e^{-s}   Dg(\mathcal{W}_s)[ \w ]ds  
             - \E \int_0^\infty \frac{e^{-2s}}{\sqrt{1- e^{-2s}}} Dg(\mathcal{W}_s)[ Z]ds.
 \end{align*}
Comparing with~\eq{eq:lhsstneq}, it only remains to show 
 \begin{align} \label{LLHS}
  \E \int_0^\infty \frac{e^{-2s}}{\sqrt{1- e^{-2s}}} Dg(\mathcal{W}_s)[ Z]ds 
             \Eq  \E \int_0^\infty e^{-2s} D^2g(\mathcal{W}_s)[Z', Z'] ds.
  \end{align}
We claim that, after swapping the order of integration, the integrands are equal. To see this, for 
fixed $\w\in\mathbb{D} {^p},$ and $s\in[0,\infty)$, we write
\[
h(\widehat{\w} )  \Eq  g(e^{-s} \w + \sqrt{1- e^{-2s} } \widehat{\w} ).
\]
Then for any $x,y\in\mathbb{D}^p$,
\ba{
D  h(\widehat{\w} )[x] & \Eq   \sqrt{1- e^{-2s} }   D  g(e^{-s} \w + \sqrt{1- e^{-2s} } \widehat{\w} )[x],\\
D^2  h(\widehat{\w} )[x, y] & \Eq   (1- e^{-2s} )   D^2  g(e^{-s} \w + \sqrt{1- e^{-2s} } \widehat{\w} )[x,y].
}
Now, the Stein equation~\eq{eq:charop1} with~\eq{eq:2ordt} implies
\[
\IE \bcls{D  h(Z )[Z] } \Eq  \E\bcls{  D^2 h (Z) [ Z', Z']   },
\]
which, using the definition of $h$, is the same as
\[
\E\big[ D  g(\mathcal{W}_s)[Z ] \big]  \Eq    \sqrt{1- e^{-2s} }   \E\big[   D^2  g(\mathcal{W}_s)[ Z',Z'] \big].
\]
This implies~\eqref{LLHS}, and thus~\eqref{SteinEq}. 
\end{proof}

\begin{remark}\label{Stein-Gaussian-method}
Theorem~\ref{prop:stnsol} can be used to establish quantitative approximation of a process~$W \in \ID^p$
by a continuous Gaussian process~$Z$ in ways typical of Stein's method.
 Taking any test function~$g \in M$, the difference $\IE\bcls{g(W)} -\IE\bcls{g(Z)}$ can be bounded by
\be{
  \Bigl|\IE\Bigl\{ \int_{\IR^d} D^2 f_g(W)[ J_u\uts]\Lambda(du)-D f_g(W)[W] \Bigr\}\Bigr|,
}
where the functions~$J_u$ and the measure~$\Lambda$ are as in the representation~\Ref{eq:covform} of the covariance 
function of~$Z$.
The quantity $\IE\{D f_g(W)[W]\}$ can then be treated in one of a number of standard ways, depending on the context.
For Theorem~\ref{thm:stnpalm}, $W$ is an integral with respect to a point process, the expectation is
evaluated using Palm theory, and~$\Lambda$ and the functions~$J_u$ emerge naturally in the resulting calculations.
\end{remark}

For making estimates when exploiting the above approach, the following two lemmas are often
useful.  They are needed, for example, in proving Theorem~\ref{thm:stnpalm}.

\begin{lemma}\label{lem:bilnbds}
If $f\in M$ and $w, w', w_1, w_2\in \mathbb{D}^p$, then
\ba{
\babs{D^2 f(w)[w_1,w_2]}& \Le  3 \norm{w_1}\norm{w_2} \norm{D^2 f(w)},
} 
and
\ba{
\babs{D^2 f(w+w')[w_1,w_2]-D^2f(w)[w_1,w_2]}& \Le  3 \norm{w_1}\norm{w_2} \norm{D^2 f(w+w') - D^2f(w)}
}
\end{lemma}

\begin{proof}
Using bilinearity, we have
\bes{
D^2 f(w)[w_1,w_2]
	& \Eq \frac{1}{2}\bclr{D^2 f(w)[w_1+w_2,w_1+w_2]-D^2 f(w)[w_1,w_1]-D^2 f(w)[w_2,w_2]}.
}
Taking the absolute value and using the triangle inequality implies  
\bes{
\babs{D^2 f(w)[w_1,w_2]}
	& \Le  \frac{1}{2} \norm{D^2 f(w)} \bclr{\norm{w_1}^2+\norm{w_2}^2+\norm{w_1+w_2}^2}\\
	& \Le  \frac{3}{2} \norm{D^2 f(w)} \big(  \norm{w_1}^2+\norm{w_2}^2 \big),
}
and we deduce from the bilinearity that for any $t>0$
\begin{align*}
 \big\vert D^2 f(w)[w_1,w_2]  \big\vert  & \Eq \babs{D^2 f(w)[t w_1, t^{-1}w_2]}  
             \Le \frac{3}{2} \norm{D^2 f(w)} \big( t^2 \norm{w_1}^2+ t^{-2}\norm{w_2}^2 \big).
\end{align*}
Taking $t^2 = \| w_2\| / \| w_1\|$ yields the first inequality.  
The second inequality follows from the same arguments, with $D^2f(w)$ replaced by  $D^2 f(w+w')-D^2f(w)$. \qedhere


\end{proof}

\begin{lemma}\label{lem:taylor}
If $f\in M$ and $w_1,w_2, J\in \mathbb{D}^p$, 
then
\bes{
   Df(w_2)[J]-&Df(w_1)[J]  \Eq D^2 f\bclr{w_1}[J,w_2-w_1] \\
	&\qquad +\int_0^1 \bbclr{D^2 f\bclr{w_1+t(w_2-w_1)}[J,w_2-w_1]-D^2 f\bclr{w_1}[J,w_2-w_1]} dt.
}
\end{lemma}

\begin{proof}
Set $h(t)=Df\bclr{w_1+t(w_2-w_1)}[J]$, and note that 
\[
h'(t) \Eq D^2 f\bclr{w_1+t(w_2-w_1)}[J,w_2-w_1]
\] is continuous on $[0,1]$. Therefore, 
\be{
h(1)-h(0) \Eq h'(0)+\int_0^1 (h'(t)-h'(0)) dt,
}
which is the lemma.
\end{proof}

We now turn to proving our second main result.

\begin{proof}[Proof of Theorem~\ref{thm:stnpalm}] Recall that
\[
    \widetilde{X}(s)  \Eq  \int_{\mathbb{R}^d} J_v(s)\, \frac{(N-\lambda)(dv)}{\sigma}\,,
\]
and let $f=f_{g}$ be the Stein solution in Theorem~\ref{prop:stnsol}.   In view of  \eqref{SteinEq}, it suffices 
to show the bound for 
\ben{\label{eq:stn2bd2}
  \babs{\IE \cA f(\tX)} \Eq \bbabs{\IE\int_{\IR^d} D^2 f(\tX)[\hJ_u\uts] \Lambda(du) - \IE\bcls{ D f(\tX)[\tX] }}. 
}
{Using  the definition of~$\tX$, we first write }
\ban{
\IE Df(\tX)[\tX] 
	& \Eq \IE \int_{\IR^d}  D f(\tX)[J_{u}]\tN(du)  \notag\\
	& \Eq \sigma^{-1}\IE \bbbcls{ \int_{\IR^d}  D f(\tX)[J_{u}]N(du)
		- \int_{\IR^d} D f(\tX)[J_{u}]\lambda(du)  }.	 \label{eq:d1exp}
}
Now, with $g(N, u):=D f(\tX)[J_{u}]$, we can write
\bes{
\IE \int_{\IR^d} g\bclr{N, u} N(du) 
	& \Eq  \IE \int_{\IR^d} g\bclr{N\sp{u}+\delta_{u}, u} \lambda(du) \\
	& \Eq  \IE \int_{\IR^d} D f\bclr{\sigma^{-1} (X\sp{u}+J_{u}-\IE[ X])}[J_{u}]\lambda(du).
}
Combining this with~\eq{eq:d1exp}, we find that
\bes{
\IE\bcls{ Df(\tX)[\tX] }
	& \Eq \sigma^{-1}\IE \bbbcls{ \int_{\IR^d} 
              \bbclr{D f\bclr{\sigma^{-1} (X\sp{u}+J_{u}-\IE[ X])}[J_{u}]-D f(\tX)[J_{u}]}\lambda(du)}.
}
Applying Lemma~\ref{lem:taylor}   with $w_1=\tX$, $w_2 =\sigma^{-1} (X\sp{u}+J_{u}-\IE[ X])$ and  $J=J_{u}$
yields
\ban{
\IE Df(\tX)[\tX]
	& \Eq \frac{1}{\sigma^2}\IE \bbbcls{ \int_{\IR^d} D^2 f(\tX)[J_{u},X\sp{u}-X+J_u] \lambda(du)}\label{eq:3ordpalm} \\
	\begin{split}
	&\quad+ \frac{1}{\sigma^2}\IE \bbbcls{\int_{\IR^d}\int_0^1 
                          \bbclr{D^2 f\bclr{\tX+t \sigma^{-1}(X\sp{u}+J_{u}-X)}[J_{u},X\sp{u}+J_{u}-X] \\
	&\qquad\hspace{4cm}-D^2 f\bclr{\tX}[J_{u},X\sp{u}+J_{u}-X]} dt\lambda(du)}. \notag
	\end{split}
}
Now the first bound of the theorem easily follows by 
 adding and subtracting
\be{
\frac{1}{\sigma^2}\IE \bbbcls{ \int_{\IR^d} D^2 f(\tX)\bclr{J_{u},\IE[X\sp{u}-X+J_u]} \lambda(du)}
}
to the right hand side of~\eq{eq:3ordpalm}, plugging the resulting expression for $\IE Df(\tX)[\tX]$ 
into~\eq{eq:stn2bd2}, 
and then applying \eqref{eq:3dstnbd}  to arrive at \eqref{eq:stnbdt3}.  

The second assertion~\eq{eq:poisstnbd} follows from the first, after observing that 
we can set \hbox{$N\sp{u}=N$,} and hence $X\sp{u}=X$.
\end{proof}




\section{${\rm M/GI/\infty}$ queue:  Proof of Theorem \ref{thm:mginfqstn}  }\label{sec:mginf}

Let us first recall some notation:  $\cS := [0,\hT]\times\mathbb{R}_+$, $\alpha$ is a finite measure on $[0,\hT]$ 
with $\a([0,\hT]) \ge 1$, and $G$ and~$\tG$ are distribution functions on $\mathbb{R}_+$. The point process~$N_n$ 
that we consider has the following form:~\eqref{eq:stnbdt2},
\[
    N_n \Eq M_n\ui + {M_n\ut}, \quad\mbox{where}\quad {M_n\ut} \Def \sum_{k=1}^{x_n} \delta_{(0,  Y_i)},
\]
where $x_n\geq 1$ is an integer, $(Y_i,\,i\geq1)$ is a sequence of \emph{i.i.d.}\ random variables with 
distribution $\tG$, and~$M_n\ui$  is a Poisson point process on~$\cS $ with intensity measure 
$\ell_n(dt,dy):=n\alpha(dt)G(dy)$ that is independent of $(Y_i,\,i\geq1)$.  
{$\Lambda$ denotes the measure $(\a \times G) + x(\d_0 \times \tG)$.}

{
In order to illustrate the use of Theorem~\ref{thm:stnpalm} in a multivariate context, we define
\[
   J_{t,y}(s) \Def \begin{cases}
                        \mathbf{1}[t\leq s<t+y]\boe\ui, &\ \mbox{if}\ t > 0;\\ 
                        \mathbf{1}[0\leq s< y]\boe\ut, &\ \mbox{if}\ t = 0,
                    \end{cases}
\]
where~$\boe\uii$, $i=1,2$, denotes the coordinate vectors in~$\IR^2$, and then define 
\[
    U_n\uii(s) \Def \int_{\cS} J_{t,y}(s)  M_n\uii(dt,dy),\quad i = 1,2;\qquad U_n \Def U_n\ui + U_n\ut.
\]
Then $X_n := (1,1)^{\tr} U_n$ models the number of customers in an ${\rm M/G/\infty}$ queue, and~$U_n$ distinguishes
those who were in the queue at time~$0$ and those who arrived afterwards.  We quantify the convergence
of $\tU_n := n^{-1/2}(U_n - \ex U_n)$ to the bivariate centred Gaussian process~$\hZ$ with covariance matrix
\[
   \ex\{\hZ(s_1)\hZ(s_2)^\tr\} \Def \int_0^{s_1} \bclr{1-G(s_2-t)} \alpha(dt)\boe\ui(\boe\ui)^\tr 
		+  x \, \tG(s_1)\bclr{1- \tG(s_2)}\boe\ut(\boe\ut)^\tr,
\]
and use this to deduce Theorem~\ref{thm:mginfqstn}; in particular, see~\Ref{ADB-Z-hat-approx}.
}

We start with the following proposition, which states  
the well known families of reduced Palm couplings  
$(M_n\sp{1,t,y}, M_n\ui)_{(t,y)\in \cS}$ and {$(M_n\sp{2,y}, M_n\ut)_{y\in\IR_+}$.}

\begin{proposition}\label{lem:mginfqcoup}
Let $M_n\uii$, $i=1,2$, be  defined  as   above. 
For $(t,y)\in \cS$, and given $M_n\ui$, let $M_n\sp{1,t,y} := M_n\ui$, $t > 0$.
Then $M_n\sp{1,t,y}$ has the reduced Palm distribution of $M_n\ui$ at~$(t,y)$.
Similarly, given~$M_n\ut$, let $Y$ be a point uniformly and independently chosen from $\{Y_1,\ldots,Y_{x_n}\}$.
Then $M_n\sp{2,y} := M_n\ut - \delta_{(0, Y)}$ has the reduced Palm distribution of~$M_n\ut$ at~$(0,y)$. 
\end{proposition}

{
Because $M_n\ui$ and~$M_n\ut$ are independent, it follows that the reduced Palm distributions of~$N_n$ are given by
$N_n\sp{t,y} \Eq N_n - \delta_{(0, Y)}\bone\{t = 0\}$, and hence that
$U_n\sp{t,y} := \int_\cS J_{t,y} N_n\sp{t,y}(dt,dy)$ satisfies
\ben{\label{ADB-reduced-Palm-U}
    U_n\sp{t,y} - U_n \Eq - J_{0,Y}\bone\{t = 0\}.
}
}


\begin{proof}[Proof of Theorem~\ref{thm:mginfqstn}]
{We start by computing the bounds in Theorem~\ref{thm:stnpalm} on the difference
$|\ex g(\tU_n) - \ex g(\hZ)|$.}  We first show that the process~$\hZ$ has a 
continuous modification.  Since, by assumption, neither $G*\alpha$ nor~$\tG$  have atoms 
in $[0,\hT]$, it is clear from~\eq{eq:minfzcov} that the covariance function is continuous. Moreover, for 
any $s\geq0$ and $0\leq h \leq (\hT-s)\wedge1$,~\eq{eq:minfzcov} and the H\"older continuity of $G*\a$
easily imply that 
\ba{
&{\bigl|\IE\bcls{\hZ(s+h)^\tr \hZ(s+h) - \hZ(s)^\tr \hZ(s+h)}\bigr|} \\
	&\quad \Eq \int_s^{s+h} \bclr{1-G(s+h-t)} \alpha(dt) + x\bclr{1-\tG(s+h)} \bcls{\tG(s+h)- \tG(s)} \\
        &\quad\Le   (c_{G,\a}+ xc_{\tG}) h^{\bo},
}
and,  similarly, that
\be{
   {\bigl|\IE\bcls{\hZ(s)^\tr \hZ(s) - \hZ(s)^\tr \hZ(s+h)}\bigr|} \Le   (c_{G,\a}+ xc_{\tG}) h^{\bo},
}
so that~\eq{eq:gpkolcont} is satisfied. Therefore, in view of   Remark~\ref{rem:checkz}, the Gaussian process~$\hZ$   
has a continuous modification. In what follows, we will work with this continuous Gaussian process, that we still 
denote by~$\hZ$.

\ignore{
Setting $\tX_n\uii(s) := n^{-1/2}(X_n\uii(s) - \ex X_n\uii(s))$ for $i=1,2$, and considering the last part of 
Proposition~\ref{lem:mginfqcoup}, along with the definitions of $\widehat J_{t,y}$ and $J_{t,y}$, and the form of 
the errors \eqref{eq:stnbdt1}--\eqref{eq:stnbdt3}, we may separately apply Theorem~\ref{thm:stnpalm} to
bound the errors between 
$\tX_n\ui$ and $Z_1 + Z_2$ as in Remark~\ref{rem:mginfgpform}, and~$\tX_n\uo$ and $Z_3$ as in Remark~\ref{rem:mginfgpform}, 
and add the contributions to obtain~\Ref{eq:mginftfbd}.
}

{
The contributions to the integrals in \eqref{eq:stnbdt1}--\eqref{eq:stnbdt3} from $\{\{0\}\times\IR_+\}$
and $\{(0,\hT] \times \IR_+\}$ can be separately bounded, and the results added for the overall bounds.
First, on $\{(0,\hT] \times \IR_+\}$,  $U_n\sp{(t,y)} - U_n = 0$, in view of~\eqref{ADB-reduced-Palm-U},
so that there is no contribution from~\eqref{eq:stnbdt2}, or from
\eqref{eq:stnbdt1} either, since $\Lambda = \a \times G = n^{-1}\ell_n$ on $\{(0,\hT] \times \IR_+\}$; and \eqref{eq:stnbdt3}  
contributes at most
\ben{\label{ADB0}
  \frac{\norm{g}_M }{ 2\sqrt{n}}   \int_{\cS} \| J_{t,y}\|^3 \alpha(dt) G( dy)  \Le  
          \frac{\norm{g}_M \alpha( [0,\hT]) }{ 2\sqrt{n}}\,. 
}
}

{
Next, on $\{\{0\} \times \IR_+\}$,  $U_n\sp{(0,y)} - U_n = -J_{0,Y}$, and so
\ignore{
Hence it only remains to use Theorem~\ref{thm:stnpalm} to show that
\ben{\label{ADB1}
     |\ex g(\tX\uo) - \ex g(Z\uo)| 
            \Le (3/2)\norm{g}_{M'}\bbbcls{\babs{n^{-1}x_n-x} + \sqrt{\pi/ (2x_n) }} 
                  +  \frac{x_n \norm{g}_M}{2n^{3/2}},
}
where $Z\uo = Z_3$ as in Remark~\ref{rem:mginfgpform}.
Now Proposition~\ref{lem:mginfqcoup} implies that 
\ben{\label{eq:condmndiffq}
    (X_n\uo, X_n\sp{0,0,y}-X_n\uo)\ \eqlaw\ (X_n\uo , - J_{0,Y}), 
}
and therefore that
}
\ben{\label{eq:mndiffq}
    \IE[U_n\sp{0,y} - U_n] \Eq - \int J_{0,y} \tG(dy) \Eq - (1- \tG)\boe\ut. 
}
Now consider the contribution to~\eqref{eq:stnbdt1}, with $\Lambda = x(\d_0 \times \tG)$ and
$\l = x_n(\d_0 \times \tG)$  on $\{\{0\} \times \IR_+\}$, taking $\hJ_{(0,y)} := J_{(0,y)} - (1-\tG)\boe\ut$.  
The contribution can be written as
\ba{
   x\int_{\IR_+} &D^2 f(\tU_n)\bcls{{\clr{J_{0,y}-(1-\tG)\boe\ut}\uts}} \tG(dy) - 
            \frac{x_n}n \int_{\IR_+} D^2 f(\tU_n)\bcls{J_{0,y},J_{0,y}-(1-\tG)\boe\ut} \tG(dy) \\
    &=\ \frac{nx - x_n}{n}\int_{\IR_+} D^2 f(\tU_n)\bcls{{\bclr{J_{0,y}-(1-\tG)\boe\ut}\uts}} \tG(dy) \\
    &   \qquad\qquad\qquad - \frac{ x_n}{n}\int_{\IR_+} D^2 f(\tU_n)\bcls{(1-\tG)\boe\ut,J_{0,y}-(1-\tG)\boe\ut} \tG(dy)\\
    &=\  \frac{nx - x_n}{n}   \int_{\IR_+} D^2 f(\tU_n)\bcls{{\clr{J_{0,y}-(1-\tG)\boe\ut}\uts}} \tG(dy),
}  
where the last line uses~\Ref{eq:mndiffq}, as well as Lemma~\ref{lem:intlin} 
below (noting in particular that \hbox{$J_{0,y}(s)= (\bone\{s\geq 0\}-\bone\{s\geq y\})\boe\ut$} 
and that~$f$ inherits from~$g$ either its smoothness property~\eq{eq:smoothgginf} or its being a function 
of a finite number of values of its argument, using  \eqref{eq4.3} of Theorem~\ref{prop:stnsol}).
}
Therefore, using~\eq{eq:2dstnbd} of Theorem~\ref{prop:stnsol}, the {contribution to~\eq{eq:stnbdt1} 
on $\{\{0\} \times \IR_+\}$} is bounded by
\be{
 (3/2)\norm{g}_{M'}\babs{n^{-1}x_n-x}.
}

For the {contribution to~\eq{eq:stnbdt2} on $\{\{0\} \times \IR_+\}$,} we use~\eq{eq:2dstnbd} and~\Ref{ADB-reduced-Palm-U},
giving
\ba{
\IE \bbbcls{\int_{\IR_+} \bbabs{D^2 &f(\tU_n)\bcls{J_{0,y}, (U_n\sp{0,y}-U_n)-\IE[U_n\sp{0,y}-U_n]} }
         n^{-1}x_n\tG(dy)}\\
	&\leq\ (3/2)\norm{g}_{M'}  
               \int_{\IR_+} \norm{J_{0,y}}
               {\IE\bbbcls{\bbnorm{\dfrac{1}{x_n}\sum_{i=1}^{x_n} \bclr{J_{0,Y_i} - \IE[J_{0,Y}]}}}n^{-1}x_n\tG(dy)}.
}
Then~\eq{eq:mndiffq} implies that
\ba{
&   {\int_{\IR_+} \norm{J_{0,y}}
               \IE\bbbcls{\bbnorm{\dfrac{1}{x_n}\sum_{i=1}^{x_n} \bclr{J_{0,Y_i} - \IE[J_{0,Y}]}}}n^{-1}x_n\tG(dy)} \\
	&=\ \frac{x_n}{n}\int_{\IR_+}\IE\bbbcls{ \bbnorm{\dfrac{1}{x_n}\sum_{i=1}^{x_n} J_{0,Y_i}-(1-\tG){\boe\ut}}}\tG(dy)  
           \Eq \frac{x_n}{n}\IE\bbbcls{\bbnorm{\dfrac{1}{x_n}\sum_{i=1}^{x_n} \mathbf{1}[Y_i  >  \cdot ] -  (1-\tG)}}   \\
	&=\  \frac{x_n}{n}\IE\bbbcls{\bbnorm{\dfrac{1}{x_n}\sum_{i=1}^{x_n} \mathbf{1}[Y_i  \leq \cdot ] - \tG}}.
}
Thus we must bound the mean of the sup-norm of the difference between an empirical CDF and its limit. 
According to 
 \cite[Corollary~1]{Massart1990} (improving on \cite{Dvoretzky1956}), for any $y>0$, 
 \be{
 \IP\bbbclr{\bbnorm{\dfrac{1}{x_n}\sum_{i=1}^{x_n} \mathbf{1}[Y_i \leq \cdot ] - \widetilde  G}> y} \Le 2 e^{-2 x_n y^2 },
 }
so that 
\be{
\IE\bbbcls{\bbnorm{\dfrac{1}{x_n}\sum_{i=1}^{x_n} \mathbf{1}[Y_i \leq \cdot ] -  \tG}} 
            \Le 2\int_{0}^\infty  e^{- 2 x_n y^2 } dy \Eq \sqrt{\frac{\pi}{2x_n}},
}
{giving a contribution to~\eq{eq:stnbdt2} on $\{\{0\} \times \IR_+\}$ of at most $(3/2)\norm{g}_{M'}\sqrt{\pi/(2x_n)}$.}

For the {contribution to~\eq{eq:stnbdt3} on $\{\{0\} \times \IR_+\}$, we use~\eq{eq:mndiffq} and~\Ref{ADB-reduced-Palm-U}}
to easily find that
\ba{
    \int_{\IR_+}\norm{J_{0,y}}\norm{U_n\sp{0,y}- U_n + J_{0,y}}^2 n^{-1}x_n \tG(dy) 
	& \Eq  \int_{\mathbb{R}_+} \| J_{0,y} \| \bigl\| J_{0,y} - J_{0,Y} \bigr\|^2 n^{-1}x_n \tG(dy) \\
         &  \Le  n^{-1} x_n,
}
{giving a contribution to~\eq{eq:stnbdt3} on $\{\{0\} \times \IR_+\}$ of at most $(1/2)\norm{g}_{M'}n^{-3/2}x_n$.
Collecting these bounds, we deduce that, for any $g \in M'$, 
\ben{\label{ADB-Z-hat-approx}
     |\ex[g(\tU)] - \ex[g(\hZ)]| \Le \norm{g}_{M'} \ppsi(x,x_n,\a,\hT),
}
where $\ppsi_n(x,x_n,\a,\hT)$ is as given in~\Ref{ADB-psi-def}.  Noting that, for $\tg\colon \ID^1 \to \IR$ 
and~$g\colon \ID^2 \to \IR$ defined by $g(w_1,w_2) := \tg(w_1+w_2)$, we have $\|\tg\|_{M'} \le 2^{3/2}\|g\|_{M'}$,
the bound~\Ref{eq:mginftfbd} in Theorem~\ref{thm:mginfqstn}  follows.}

\medskip

To prove the bound on the L\'evy--Prokhorov distance, 
we use the main results of \cite{Barbour2021}, as stated in 
Theorem~\ref{thm:smooth} below. The first hypothesis of the theorem is satisfied with $\kappa_2=0$, and with $\kappa_1$ 
upper bounded by a quantity of order  $\bigo(|n^{-1}x_n - x| + \a([0,T]) n^{-1/2})$, read from the 
bound~\eq{eq:mginftfbd} just established (noting that $\norm{g}_{M'}\leq \norm{g}_{M^0}$). 

To bound the modulus of continuity terms, we use Lemma~\ref{lem:modc}, {treating the components $\tX_n\ui$ and~$\tX_n\ut$
of~$\tU_n$ separately (so that $\tU_n\uii = \tX_n\uii \boe\uii$).} To verify~\eq{eq:modch1} for~$\tX_n\ui$, let
 \besn{\label{eq:regions}
 \cR_s &\Def \bclc{(u,y): 0 \le u \leq s,\ 0 < y < s-u}, \quad \cR_1(s_1,s_2) \Def \cR_{s_2}\setminus \cR_{s_1} \\
     &\qquad\mbox{and}\quad \cR_2(s_1,s_2) \Def (s_1,s_2]\times \IR_+,
      \quad s_1 < s_2.
 }
Fix $0\leq s < t \leq \hT$ with $1/(2n) \leq (t-s)\leq  1/2$.
Recalling the definition of the random measure $M_n\ui$, we have 
\be{
  \tX_n\ui(t)-\tX_n\ui(s) \Eq  \tY_n(2;s,t) - \tY_n(1;s,t) ,
}
where 
\begin{align}\label{YST}
   \tY_n(i;s,t) \Def n^{-1/2}\{M_n\ui(\cR_i(s, t))  - \ell_n(\cR_i(s,t))\}, \,\, i=1,2,
\end{align}
are (dependent) centred and normalized Poisson random variables with means 
\ben{\label{ADB3}
    \ell_n(\cR_1(s,t)) \Le nc_{G,\a}(t-s)^\b \quad\mbox{and}\quad \ell_n(\cR_2(s,t)) \Le n c_\a (t-s)^\b,
}
by the H\"older continuity of $G*\a$ and~$A$.  Now, for~$W_n$ a sum of $n$ independent Bernoulli
random variables with success probability $\tilde p \leq \mu$,
it follows from Rosenthal's inequality that, for any $l \ge 1$, 
\ben{\label{ADB3.5}
      n^{-l}\ex|W_n-\IE[W_n] |^{2l} \Le C_{2l}n^{-l}  \max\{(n\m)^{l} , n\m\} \Le C_{2l} \max\{\m^{l}, n^{-l+1}\mu\},
}
where~$C_{r}$ is the Rosenthal constant for exponent~$r$. 
 {A limiting argument shows that the inequality~\eqref{ADB3.5} holds also for $W_n \sim \Po(n\m)$.}
Thus 
it follows that, for any $l \ge 1$ and $|t-s| \ge (1/2) n^{-1/\b}$,
\begin{align}
      \pr\big[|\tY_n(i;s,t)| \geq \th/2 \big]& \Le C_{2l}2^{2l}\th^{-2l}\max\big\{(c_{G,\a}\vee c_\a)^l  |t-s|^{l\b}, (c_{G,\a}\vee c_\a) n^{-l+1}  |t-s|^{\b}\big\} \notag \\
	&               \Le K_l\ui \th^{-2l}|t-s|^{l\b},\quad i=1,2, \label{in_view51}
\end{align}
which implies~\eq{eq:modch1} for~$\tX_n\ui$ with $\MM = n^{1/\b}$, $ K = 2 K_l\ui$, $a = l\b-1$ and $b=2l$, 
for any $l \ge 1$.
\ignore{
 If $0 < \ell_n(\cR_i(s,t)) \le 1$, it follows that 
\[
    \pr[|\tY_n(i;s,t)| \geq \th/2] \Le 2C_{2l}2^{b_l}\th^{-2l} n^{-l} \Le 2^{1+l+b_l}C_{2l}\th^{-2l}|t-s|^l,   
\]
this last because $n(t-s) \ge 1/2$.  Since $\b \le 1$, this establishes~\eq{eq:modch1} for~$\tX_n\ui$,
for any $l \ge 1$,
with $a_l$ and~$b_l$ for $a$ and~$b$, and with $K = K_l\ui := 2C_{2l}\max\{(c_{G,\a}\vee c_\a)^l,2^l\}$.
}

To establish~\eq{eq:modch2} for~$\tX_n\ui$, with $\MM = n^{1/\b}$, observe that, for all~$s$ such that 
$(k-1)/\MM \le s \le k/\MM$,
\ba{
     \Big\vert \tX_n\ui(s) &- \tX_n\ui \bclr{(k-1)/\MM} \Big\vert  \\
                 &\Le n^{-1/2}\sum_{i=1}^2\Big\{M_n\ui \Bigl( \cR_i\Bigl(  \frac{k-1}{\MM},\frac{k}{\MM} \Bigr) \Bigr) 
                                               + \ell_n\Bigl( \cR_i\Bigl(  \frac{k-1}{\MM},\frac{k}{\MM} \Bigr) \Bigr) \Big\} \\
                 &\Eq   \sum_{i=1}^2\tY_n\Big(i; \frac{k-1}{\MM},  \frac{k}{\MM}\Bigr)  
                                    + 2n^{-1/2}\sum_{i=1}^2\ell_n\Big(  \cR_i\Bigl( \frac{k-1}{\MM},\frac{k}{\MM}\Bigr) \Big),
}
and that $\ell_n(\cR_i( \frac{k-1}{\MM}, \frac{k}{\MM})) \Le  c_{G,\a}\vee c_\a$, see \Ref{ADB3} and \Ref{YST}. 
Hence, if  $\theta \geq 4(c_{G,\a}\vee c_\a)n^{-1/2}$, then
\begin{align*}
    \pr\left[\sup_{ \frac{k-1}{\MM} \le s \le \frac{k}{\MM} }|\tX_n\ui(s) - \tX_n\ui \bclr{(k-1)/\MM}| \ge \th\right]
         & \Le \pr\left(    \sum_{i=1}^2 \Bigl\vert \tY_n\Bigl(i; \frac{k-1}{\MM},  \frac{k}{\MM}\Bigr) \Bigr\vert   \ge  \th/2 \right) \\
         &\Le  \sum_{i=1}^2 
               \pr\biggl[  \Bigl\vert \tY_n\Bigl(i; \frac{k-1}{\MM},  \frac{k}{\MM}\Bigr) \Bigr\vert   \ge  \th/4 \biggr],
\end{align*}
and this probability is bounded by 
$2^{2l'+1}K_{l'}\ui\th^{-2l'}n^{-l'}$, for any $l' \ge 1$,
in view of  \eqref{in_view51}  as established above.
Hence, for $\th \ge 4(c_{G,\a}\vee c_\a)  n^{-1/2}$, we can take
\ben{\label{ADB4}
      \f_\MM\ui(\th) \Def n^{ \frac{1}{\b} -l'}2^{2l'+1}K_{l'}\ui\th^{-2l'} \Le 2^{2l'+1}K_{l'}\ui\th^{-2l'}\ee^{l'\b-1}
}
in~\eq{eq:modch2}, for any $l' \ge 1/\b$ and $\ee \ge n^{-1/\b}$, to be compared with the bound in~\eq{eq:modch1}.  
In particular, taking $\th \ge 4(c_{G,\a}\vee c_\a)n^{-1/2}$  and $l' = l$, for any $l > 1/\b$,
and applying Lemma~\ref{lem:modc} with $\MM = n^{1/\b}$, 
it follows that, for any $\ee \in (n^{-1/\b},1]$ and for any $\th > 4(c_{G,\a}\vee c_\a)n^{-1/2}$, we have
\ben{\label{ADB4.5}
   \IP[\omega_{\tX_n\ui} (\epsilon) \geq \theta/2] 
         \Le    \hT C
             \theta^{-{2l}} \epsilon^{l\beta-1},
}
for a suitable constant~$C$ that does not depend on $(\epsilon, n, \theta)$.  {By observing that, for $\ee \ge n^{-1/\b}$
and $\th \leq 4(c_{G,\a}\vee c_\a)n^{-1/2}$, the bound \eqref{ADB4.5} is comparable to or larger than~$1$, the constant~$C$
can be chosen in such a way that the bound is valid for all~$\th > 0$.}

 Turning to $\tX_n\uo$,  we can first assume  that $x_n/n \leq 2 x$  without loss of generality, because $\hT \ge n^{1/2}|x_nn^{-1} - x|$ and the final bound is only meaningful for $\hT\ll n^{1/8}$. Note that, for $s < t$, 
\[
   \tX_n\uo(t) - \tX_n\uo(s) \Eq n^{-1/2}\sum_{i=1}^{x_n}\{ \mathbf{1}[s < Y_i \le t] - (\tG(t) - \tG(s))\}
\]
is a normalized sum of independent centred Bernoulli random variables, and that, by assumption,
$x_n (\tG(t) - \tG(s)) \le 2n x c_{\tG}(t-s)^\tbe$.  Arguing exactly as for $\tX\ui$ now yields~\eq{eq:modch1} for~$\tX_n\uo$,
for any $l \ge 1$, 
with $1 + a = l\tbe$ and $b = 2l$, and with $K = K_l\uo = 2 ^{l+1} x C_{2l} c_{\tG}\max\{x c_{\tG},1\}^{l-1}$.
For~\eq{eq:modch2}, the argument is again as for~$\tX_n\ui$.   We first write 
\ba{
 & \sup\left\{   |\tX_n\uo(s) - \tX_n\uo \bclr{(k-1)/\MM}| :   \frac{k-1}{M} \leq s \leq \frac{k}{M} \right\} \\
            & \qquad\qquad \Le n^{-1/2}\sum_{i=1}^{x_n}\Bigl\{ \mathbf{1}\Bigl[ \frac{k-1}{\MM} < Y_i \le \frac{k}{\MM} \Bigr] 
              + \Bigl[  \tG\Bigl(\frac{k}{\MM}\Bigr) - \tG\Bigl( \frac{k-1}{\MM} \Bigr) \Bigr]  \Bigr\} \\
             &\qquad\qquad \Le  n^{-1/2} \big\vert W_n - \IE[W_n] \big\vert  +  4x n^{1/2}\{\tG(k/\MM) - \tG((k-1)/\MM)\},
}
 where $W_n$ is a sum of $x_n$ i.i.d.\ Bernoulli random variables with success probability 
 $\tilde{p} \leq c_{\widetilde{G}} M^{-\beta}$.
Therefore, with $M=n^{1/\beta}$, it follows from \eqref{ADB3.5} and by first considering $\th >  {8x}c_{\tG}n^{-1/2}$ that 
we can take 
\ben{\label{ADB5}
               \f_\MM\uo(\th) \Def n^{-l' + \frac{1}{\b}}2^{2l'}K_{l'}\uo\th^{-2l'}
}
for $\f_\MM(\th)$ in~\eq{eq:modch2}  for any $l' \ge 1$.
Hence, from Lemma~\ref{lem:modc}, 
 for any $\ee \in (n^{-1/\b},1]$
and any $\th > 0$,  we have
\ben{\label{ADB5.5}
   \IP[\omega_{\tX_n\uo} (\epsilon) \geq \theta/2] 
         \Le    \hT C  \theta^{-{2l}} \epsilon^{l\b-1},
}
for any $l \ge 1/\b$, for a suitable constant $C$.

 For the analogous inequality for~$\hZ$,
an easy calculation shows that, for any $0\leq u < s \leq \hT$, 
\be{
  \IE\bcls{ \bigl\vert  {\hZ(s) - \hZ(u)} \bigr\vert ^2} \Le 2(1+x) c(s-u)^\bo,
}
and so \cite[Remark~1.6]{Barbour2021}
implies  there is a constant $C$ depending on $x, c, \bo$ such that, for each component $\hZ\uii$ of~$\hZ$,
\ben{\label{modczmg}
   {\IP[\omega_{\hZ\uii} (\epsilon) \geq \theta]  \Le  C\hT \theta^{-{2l}} \epsilon^{l\bo-1},\quad i = 1,2,}
}
for any $l \ge 1$.

{A bound on the L\'evy--Prokhorov distance between $\law(\tU_n)$ and~$\law(\hZ)$}
now follows by  using \Ref{ADB4.5}, \Ref{ADB5.5} and~\eq{modczmg} in Theorem~\ref{thm:smooth} below, 
with $\kappa_1 = \bigo\bclr{|x_n n^{-1}-x| + \alpha\bclr{[0,\hT]}n^{-1/2}}$ {(from~\Ref{ADB-Z-hat-approx})}
and $\kappa_2=0$, giving a bound of order
\begin{align}\label{O_1}
   \bigo\Bigl( \th + \d\sqrt{\hT\log n} + (\ee\d)^{-3}\k_1\hT^{3/2}
                    + \hT\theta^{-{2l}} \ee^{l\bo-1} \Bigr)
\end{align}
for the L\'evy--Prokhorov distance, for any $\d,  \th > 0$ and   $\ee \in( n^{-1/\b},1)$.
 Taking $\delta \sqrt{T} = \th$ and matching $(\ee\d)^{-3}\k_1\hT^{3/2}
                    = \hT\theta^{-{2l}} \ee^{l\bo-1} $ reduces the bound \eqref{O_1} to 
\[
    \bigo\Bigl( \th \sqrt{\log n} +  \kappa_1T^3 (\ee \th)^{-3}  \Bigr)  \quad  
              \text{with}\quad \ee = \big(\kappa_1T^2 \theta^{2l-3 } \big)^{\frac{1}{2+l\b}} > n^{-1/\b}
\] 
and then balancing $\theta$ with $ \kappa_1T^3 (\ee \th)^{-3}$ yields the bound $  \bigo\bigl( \th \sqrt{\log n}   \bigr)  $ 
with 
\ba{
  & \th \Eq (\k_1\hT^{3})^{(l\bo-1)/(6l + 4l\bo - 1)}\hT^{ 3/(6l+4l\bo-1)}. 
}
That is,  we have a bound of order
\[
    \bigo \Big(\sqrt{\log n} (\hT^{4} n^{-\ff})^{(l\bo-1)/(6l + 4l\bo - 1)}\hT^{{3/(6l+4l\bo-1)}} \Big),
\]
for any $l \ge 1/\b$, where 
\[
   \ff \Def \min\left\{ \frac{1}{2} + \log_n\Bigl(\frac{\hT}{\a\bclr{[0,\hT]} }\Bigr),\, 
                                        \log_n \Bigl(\frac{\hT}{|x_nn^{-1} - x|}\Bigr) \right\}.
\]
The simplified bound given in the statement of Theorem~\ref{thm:mginfqstn}, for $\bo = 1$ and 
$n^{1/2}|x_n/n-x| \le \a\bclr{[0,\hT]} = \bigo(\hT)$,  follows, for any $\chi > 0$, by taking~$l$ large enough. \qedhere

  \end{proof}


\section{${\rm GI/GI/\infty}$ queue: Proof of Theorem \ref{thm:gginfqstn} }\label{sec:gginf}

Let us first recall some notation from Section~\ref{sec:mggres}.  The stationary renewal process~$V_n$ is driven by~$R/n$, 
and the point process~$N_n$ that we consider has the form
\[
    N_n \Def \sum_{i=1}^{\floor{n\hT}} M_n(i/n) \delta_{(i/n, Y_i)},
\]
where 
\begin{description}
\item[(i)] $M_n(i/n)$ marks the arrival of a customer at time $i/n$, and can be represented as 
\[
       M_n(i/n) \Def \sum_{j=1}^\infty \mathbf{1}[ R_0 + R_1+... + R_j = i],  
\]
where $(R_j,\, j\geq 1)$ are independent copies of~$R$, and~$R_0$ has the delay distribution 
\[
    \IP(R_0 = k) \Eq  m^{-1} \IP(R \geq k)\,, \,\,\, k\in\IN;
\]
\item [(ii)] the service times  $Y := (Y_i,\,i\geq 1)$ are i.i.d.\ with distribution~$G$, 
\item [(iii)] $M_n$ and~$Y$ are independent.
\end{description}

As in the previous section, we begin with a coupling lemma.

\begin{lemma}\label{lem:gginfpalm}
With the above notation,  let $(R_{ik},\,1 \le i,k <\infty)$ and $(R_{ik}',\,1 \le i,k <\infty)$ be independent 
i.i.d.\ sequences with the same distribution as~$R$ that are also independent of $V_n$ and~$Y$.  
Define  
\ba{ 
   S_{ij} &\Def \sum_{k=1}^j R_{ik} \quad  \mbox{and} \quad   S_{ij}' \Def\sum_{k=1}^j R_{ik}',
    \qquad i,j\in\{1,2,\ldots\}; \\
   \hM_n\sp{i/n} &\Def \sum_{j=1}^\infty\bclr{ \delta_{(i+S_{ij})/n} + \delta_{(i-S_{ij}')/n} },
    \qquad i\in\{0,1,...,n\};
}
then set 
\ba{
T_{i} &\Def \min\bbclc{\floor{n\hT}-i+1,\, 
            \inf\bclc{j\geq 1\colon M_n\bclr{(i+j)/n}= \hM_n\sp{i/n}\bclr{(i+j)/n} = 1}}, \\
T_{i}' &\Def \min\bbclc{i,\, \inf\bclc{ j\geq 1\colon M_n\bclr{(i-j)/n}=\hM_n\sp{i/n}\bclr{(i-j)/n} = 1}};
}
and finally define
\be{
   M_n\sp{i/n}(j/n) \Def
        \begin{cases}
          \hM_n\sp{i/n}(j/n), & i- T_{i}' <j < i+  T_{i} \\
           M_n(j/n), & \mbox{otherwise}.
        \end{cases}
}
Then
\be{
    N_n\sp{i/n } \Def \sum_{j=1}^{\floor{n\hT}} M_n\sp{i/n}(j/n) \delta_{(j/n, Y_j)}
}
has the reduced Palm distribution of~$N_n$ at $(i/n,y)\in\cS$. 
\end{lemma}

\begin{proof}
It is well known that $\hM_n\sp{i/n}$ has the reduced Palm measure of $M_n$ at $(i/n)$; see for 
example \cite[Chapter~13]{Daley2008}.  Then note that $M_n\sp{i/n}\eqlaw \hM_n\sp{i/n}$, since~$T_i$ is the 
first time there is a renewal in both the zero delayed renewal process started from~$i$: $(i+S_{ij})_{j\geq 1}$ 
appearing in the definition of $\hM_n\sp{i/n}$, and the analogous stationary process induced by $M_n$, at which 
point we can continue using either process without changing the distribution. A similar statement holds for $T_i'$ 
but now moving backwards in time. From this observation, it is clear that $N_n\sp{i/n}$ is distributed as claimed. 
\end{proof}

\begin{remark}
We write $N_n\sp{i/n}$ for the reduced Palm distribution at $(i/n,y)$, for all~$y$.
\end{remark}

The next result gives an expression for $\IE[X_n\sp{i/n}-X_n]$, to be used in~\eq{eq:stnbdt1} of 
Theorem~\ref{thm:stnpalm}.

\begin{lemma}\label{lem:gginfdiffmean}
With the notation above, let
\be{
   X_n\sp{i/n} \Def \int_{\cS}   J_{t,y}(s) N_n\sp{i/n}(dt,dy).
} 
If  $\IE[R^r]<\infty$ for some $r > 3$,  and if~$G$ satisfies the assumptions~\Ref{eq:gcon}, 
\ignore{
 functions~$(J'_{i})_{i=1}^{\floor{n\hT}}$, constants $(B_{ij})_{1\leq i, j\leq {\floor{n\hT}}}$ and a constant $C$ 
depending only on $R$ such that $\norm{J'_i}\leq C$, $\abs{B_{ij}}\leq C n^{\eps}$ and for any $\eps\in(0,1)$, 
\begin{align}
  &\IE\bcls{X_n\sp{i/n}(s)} - \IE\bcls{X_n(s)} - \bone[s\geq i/n]\overline G(s-\tsfrac{i}{n})\bbbclr{\frac{v^2-m^2}{m^2}} 
                                  \label{eq:111}\\
	&=  J_i'(s)\bclc{n^{\eps(2-r)} +n^\eps\bclr{\Delta_{n^{\eps-1}}(G)+\bone[i \leq  n^\eps]}}
  +\mathbf{1}[i> n^\eps]  \sum_{j=\floor{i-n^\eps}+1}^{\floor{i+n^\eps}}B_{ij}\bclr{ I_{\frac{j}{n}}(s)-I_{\frac{j+1}{n}}(s)},
             \notag
\end{align}
where 
\ben{\label{ADB-Delta-def}
    \Delta_h(G) \Def \sup_{0\leq x \leq \hT-h }\bclr{G(x+h)- G(x)}.
}
}
then there is a constant $C$,  depending only on $\law(R)$, $\b$ and~$g_G(0)$, 
such that
\[
     A_{n,i}(s) \Def \IE\bcls{X_n\sp{i/n}(s)} - \IE\bcls{X_n(s)} 
                   - \bone{\{s \geq in^{-1}\}}\oG(s- in^{-1})\bbbclr{\frac{v^2-m^2}{m^2}} 
\]
satisfies
\[
        |A_{n,i}(s)| \Le   C   \Bigl\{n^{-\b} + (|i-\lfloor ns \rfloor| + 1)^{-(r-2)}
                                           + {i}^{-(r-2)} \Bigr\}.
\]

 \end{lemma}

\begin{proof} First note that 
 
\ba{
    X_n\sp{i/n}(s) - X_n(s)&\Eq \sum_{j=1}^{\floor{n\hT}} (M_n\sp{i/n}(j/n) - M_n(j/n)) \mathbf{1}[j/n \leq s < Y_j + j/n],
}
so that, using the independence of $Y$ and $(M_n,M_n\sp{i/n})$, we have
\ban{
   \IE\bcls{X_n\sp{i/n}(s)} - \IE\bcls{X_n(s)}
       &\Eq \sum_{j=1}^{\floor{ns}}\Big(\IE\bcls{ M_n\sp{i/n}(j/n)} - \IE\bcls{ M_n(j/n)} \Big) \oG(s-jn^{-1}).
     \label{GGI-ADB1}
}
By stationarity,  $\IE[M_n(j/n)]=1/m$; and  by Lemma~\ref{lem:gginfpalm}, 
\[
     \IE[M_n\sp{i/n}(j/n)] \Eq u^0_{|j-i|} \Def u_{|j-i|}\bone{\{j \neq i\}},
\]
where $u_l$ is the probability that there is a renewal at time~$l$ in a renewal process with inter-arrival 
distribution~$\law(R)$, started from zero.  Hence we can write~\Ref{GGI-ADB1} in the form
\ban{
   \IE\bcls{X_n\sp{i/n}(s)} - \IE\bcls{X_n(s)} 
     &\Eq \sum_{j=1}^{\floor{ns}}(u^0_{|j-i|} - m^{-1}) \oG(s-jn^{-1}) \non\\
     &\Eq \sum_{j=1}^{\floor{ns}}(u^0_{|j-i|} - m^{-1}) \oG(s-in^{-1})\bone{\{s \geq in^{-1}\}} \label{GGI-ADB2} \\
     &\mbox{}\quad        + \sum_{j=1}^{\floor{ns}}(u^0_{|j-i|} - m^{-1})
                   \Big(   \oG(s-jn^{-1}) - \oG(s-in^{-1})\bone{\{s \geq in^{-1}\}}  \Big). \non
}
The remainder of the proof consists of showing that the first term in~\Ref{GGI-ADB2} is close to
$\bone{\{s \geq in^{-1}\}}\oG(s-in^{-1})(v^2-m^2)/m^2$, and that the second term is small.  The main tool is
the inequality 
\begin{align}\label{eqsee2}
      \babs{\IE[M_n\sp{i/n}(j/n)]-1/m} \Le C_R (\abs{j-i}+1)^{-(r-1)},
\end{align}
for a suitable constant~$C_R$, which follows from  \cite[Corollary~(6.21)]{Pitman1974}, together with the
observation that
\ben{\label{GGI-ADB3}
    \sum_{j=-\infty}^{\infty}(u^0_{|j-i|} - m^{-1}) \Eq \frac{v^2-m^2}{m^2}\,.
}

For the first term in~\Ref{GGI-ADB2}, using \Ref{eqsee2} and~\Ref{GGI-ADB3}, for $s \ge i/n$, we have
\ban{
   \Bigl| \sum_{j=1}^{\floor{ns}}(u^0_{|j-i|} - m^{-1}) - \frac{v^2-m^2}{m^2} \Bigr|
      &\Le C_R\Bigl\{ \sum_{j=-\infty}^0 (|j-i|+1)^{-(r-1)} 
                       + \sum_{j = \floor{ns}+1}^{\infty}(|j-i|+1)^{-(r-1)} \Bigr\} \non\\
      &\Le \frac{C_R}{r-2} \Bigl\{  {i}^{-(r-2)} + |\floor{ns} - i + 1|^{-(r-2)} \Bigr\}; 
                \label{GGI-ADB4}
}
for $s < i/n$, the term is zero because of the factor $\bone{\{s \geq in^{-1}\}}$.  For the second term 
in~\Ref{GGI-ADB2}, for $s \ge i/n$, using~\Ref{eq:gcon},
\[
   |\oG(s - i/n) - \oG(s - j/n)| \Le g_G(0)|n^{-1}(i-j)|^\b, \quad 1 \le j \le \floor{ns},
\]
so that
\ban{
   \Big\vert  &\sum_{j=1}^{\floor{ns}}(u^0_{|j-i|} - m^{-1})
                   \Bigl(\oG(s-jn^{-1}) - \oG(s-in^{-1})\bone{\{s \geq in^{-1}\}}  \Bigr)\Big\vert \non\\
          &\Le \sum_{j=1}^{\floor{ns}} C_R (|j-i| + 1)^{-(r-1)}\,g_G(0) n^{-\b}|i-j|^\b 
          \Le n^{-\b}C_R g_G(0) \Bigl\{{ \frac2{r-2-\b}}\Bigr\}. \label{GGI-ADB5}
}
Finally, for $s < i/n$, we have
\besn{
    \sum_{j=1}^{\floor{ns}}|u^0_{|j-i|} - m^{-1}|\, \oG(s-jn^{-1}) &\Le C_R \sum_{j=1}^{\floor{ns}}(|j-i|+1)^{-(r-1)} \\
            & \Le C_R {\Bigl\{\frac{r-1}{r-2}\Bigr\}}\,|\floor{ns} - i + 1|^{-(r-2)}. \label{GGI-ADB6}
}
Combining \Ref{GGI-ADB4}--\Ref{GGI-ADB6} with~\Ref{GGI-ADB2} proves the lemma.  \qedhere

\end{proof}

\medskip
We now use the coupling of Lemma~\ref{lem:gginfpalm} and Theorem~\ref{thm:stnpalm}
to prove Theorem~\ref{thm:gginfqstn}.

\begin{proof}[Proof of Theorem~\ref{thm:gginfqstn}]
Following  the notation in Theorem~\ref{thm:stnpalm},  we set 
\begin{center}
     $\Lambda(dt,dy) \Def dt\, G(dy)$ \qquad and \qquad 
               $\hJ_{t,y} \Def \dfrac{m}{v}\,J_{t,y}-\dfrac{m+v}{v}\,\oG_t I_{t}$,
\end{center}
 where $\oG_t(s) := \oG(s-t)$ and $I_t(s) := \bone\{s\geq t\}$. 

From the assumptions~\Ref{eq:gcon} on~$G$, it is clear from~\eq{eq:covgginfz} that the covariance 
function is continuous. 
Moreover, for any $s\geq0$ and $0\leq h \leq (\hT-s)\wedge 1$, they easily imply, with~\eq{eq:covgginfz},
that there is a constant $c'$ such that
\ba{
    \IE\bcls{\bclr{Z(s+h)-Z(s)}^2} \Le c' h^\beta,
}
so that~\eq{eq:gpkolcont} is satisfied. Therefore, in view of  Remark~\ref{rem:checkz},  we can assume that~$Z$   
has continuous sample paths. 

To bound the term corresponding to~\eq{eq:stnbdt1}, there are two issues. The first is that we need to compare 
integration against the atom-less $\Lambda$ to the atoms of $\lambda_n/\sigma_n^2$ corresponding to the 
renewals of~$X_n$ occurring  on a discrete lattice. The second is that there are error terms in the differences 
of the means of $X_n\sp{i/n}$ and~$X_n$, as given in Lemma~\ref{lem:gginfdiffmean}. To handle the first issue, 
we introduce $\Lambda_n := \frac{1}{n}\sum_{i=1}^{\floor{n\hT}} \bclr{ \delta_{i/n} \times G}$ as a discretized 
version of $\Lambda$. Then we can compute
\ban{
\IE& \bbbcls{\int D^2 f(\tX_n)[\hJ_{t,y}\uts] \Lambda(dt,dy) 
              -\int D^2 f(\tX_n)[\hJ_{t,y}\uts] \Lambda_n(dt,dy)}  \label{eq:intttgg} \\
    &=\ \sum_{i=1}^{\floor{n\hT}}\int_{\IR_+} \IE\bbbcls{\int_{(i-1)/n}^{i/n}  
      (D^2 f(\tX_n)[\hJ_{t,y}\uts]  -  D^2 f(\tX_n)[\hJ_{i/n,y}\uts])\, dt  } G(dy) 
                                      \notag \\
    &=\ \sum_{i=1}^{\floor{n\hT}}\int_{\IR_+} \IE\bbbcls{\int_{(i-1)/n}^{i/n}  
      (D^2 f(\tX_n)[\hJ_{t,y}-\hJ_{i/n,y}, \hJ_{t,y}]  +
                      D^2 f(\tX_n)[\hJ_{i/n,y}, \hJ_{t,y}-\hJ_{i/n,y}])\, dt  } G(dy).\notag
}
We work on
\ben{\label{eq:ggbdsm}
      \babs{D^2 f(\tX_n)[\hJ_{t,y}-\hJ_{i/n,y}, \hJ_{s,y}] },
}
where $t\in((i-1)/n,i/n]$, and $s\geq0$.
First, note that we can write
\ben{\label{eq:hatjrep}
    \hJ_{r,y} \Eq \frac{m}{v}\bclr{I_r- I_{r+y}} -\frac{m+v}{v}\bclr{I_r- G_r},
}
so that 
\be{
\hJ_{t,y}-\hJ_{i/n,y} \Eq 
      \frac{m}{v}\bclr{(I_t-I_{i/n})- (I_{t+y}-I_{i/n+y})} -\frac{m+v}{v}\bclr{(I_t-I_{i/n})- (G_{t} - G_{i/n})},
}
and thus bilinearity implies that \eqref{eq:ggbdsm} is bounded by
\besn{\label{eq:n2bd}
C \bbclc{\babs{D^2 f(\tX_n)&[I_t-I_{i/n}, \hJ_{s,y}] } \\
	&\mbox{}+\babs{D^2 f(\tX_n)[I_{t+y}-I_{i/n+y}, \hJ_{s,y}] }+\babs{D^2 f(\tX_n)[G_t-G_{i/n}, \hJ_{s,y}] }};
}
here and below we allow $C$  to change from line to line, but only depending on~$\law(R)$ and  $g_G(0)$.

To bound the last term of~\eq{eq:n2bd}, Equation~\eq{eq:2dstnbd} of Theorem~\ref{prop:stnsol} and the 
assumption~\eq{eq:gcon} on~$G$ (noting that $\abs{t-i/n}\leq 1/n$) imply that
\ben{\label{eq:n2bdt3}
    \babs{D^2 f(\tX_n)[G_t-G_{i/n}, \hJ_{s,y}] } \Le C \norm{g}_{M'}  n^{-\beta}.
}
To bound the first two terms, we apply \eq{eq4.3} and, noting that $J_{s,y}=I_s-I_{s+y}$, we have
\ba{
   \babs{D^2& f_g(w)[I_t-I_{i/n},J_{s,y}]} \Le \int_{0}^\infty e^{-2z} 
              \IE\bbcls{\babs{D^2 g(w e^{-z}+ \sqrt{1-e^{-2z}} Z) [I_t-I_{i/n},J_{s,y}]}} dz \\
	&  \Le
  \begin{cases}
     S_g\,n^{-1/2},  &\text{if}\ g \mbox{ satisfies~\eq{eq:smoothgginf}}, \\
     k^2 \norm{g}_{M'} \sum_{j=1}^k \bone\{t_j \in \big((i-1)/n,i/n\big]\}, 
                     &\text{if}\  g(w) = F\bclr{w(t_1),\ldots, w(t_k)},
  \end{cases}
}
where, in the first case, we use the smoothness condition~\eq{eq:smoothgginf}, and in the  second the explicit 
expression
\be{
     D^2g(w)[w_1,w_2] \Eq \sum_{j,\ell=1}^k F_{j \ell}\bclr{w(t_1),\ldots, w(t_k)} w_1(t_j) w_2(t_\ell), 
}
where we write $F_{j \ell}$ for the mixed partial derivative of~$F$ in the coordinates $j$ and~$\ell$.
Similarly, using Lemma~\ref{lem:intlin}, we can write
\ba{
\babs{D^2 f(\tX_n)&[I_t-I_{i/n},  G_{s}] } 
	\Eq \bbabs{\int D^2 f(\tX_n)[I_t-I_{i/n},  I_{s+y'}]  G( dy' )} \\
    &  \Le
     \begin{cases}
      S_gn^{-1/2},  &\text{if}\ g \mbox{ satisfies~\eq{eq:smoothgginf}}, \\
      k^2 \norm{g}_{M'} \sum_{j=1}^k \bone\{t_j \in \big((i-1)/n,i/n\big]\}, 
                    &\text{if}\  g(w) = F\bclr{w(t_1),\ldots, w(t_k)},
\end{cases}
}
and there are analogous bounds for the last two displays, when replacing $t$ by $t+y$ and $i/n$ by $i/n +y$. 
Noting that $\sum_{i=1}^{\floor{n\hT}}\sum_{j=1}^k \bone\{t_j \in \big((i-1)/n,i/n\big]\} =k$, we can apply these 
last inequalities with the representation~\eq{eq:hatjrep} for $\hJ_{s,y}$ to see that the absolute value 
of~\eq{eq:intttgg} is bounded by
\besn{\label{eq:inttggbd1}
\bbbabs{\IE\bbbcls{\int& D^2 f(\tX_n)[\hJ_{t,y}\uts] \Lambda(dt,dy) -
                         \int D^2 f(\tX_n)[\hJ_{t,y}\uts] \Lambda_n(dt,dy)}} \\
	&\leq C \hT
   \begin{cases}
      \norm{g}_{M'} n^{-\beta} + S_g n^{-1/2}, 
	&\text{if}\ g \mbox{ satisfies~\eqref{eq:smoothgginf}}, \\ 
	\norm{g}_{M'}\bclr{ n^{-\beta} + k^3 n^{-1}}, &\text{if}\  g(w) = F\bclr{w(t_1),\ldots, w(t_k)}.
   \end{cases}
}

To finish bounding~\eq{eq:stnbdt1}, we can apply Lemma~\ref{lem:intlin} (after rewriting 
$\hJ$ as per~\eq{eq:hatjrep}) to find
\ban{
\IE \bbbcls{\int D^2 f(\tX_n)[\hJ_{t,y}\uts] \Lambda_n(dt,dy)} 
	&\Eq \frac{m^2}{v^2 n}\sum_{i=1}^{\floor{n\hT}}\int_{\IR_+}
                   \IE \bbbcls{ D^2f(\tX_n)\bbcls{J_{i/n,y}\uts}}G(dy) \non\\
	&\mbox{}\qquad + \frac{v^2-m^2}{v^2 n}\sum_{i=1}^{\floor{n\hT}}
                \IE \bbbcls{ D^2f(\tX_n)\bbcls{(I_{i/n} \oG_{i/n})\uts}}. \label{GGI-ADB7}
}
Since $\lambda_n/\sigma_n^2= (m^2/v^2)\frac{1}{n}\sum_{i=1}^{\floor{n\hT}} \bclr{ \delta_{i/n} \times G}$, we have
\ban{
\IE \bbbcls{ \int D^2 &f(\tX_n)\bcls{J_{t,y},\IE[X_n\sp{t,y}-X_n+J_{t,y}]}\bclr{\sigma_n^{-2}\lambda_n(dt,dy)}} \non\\
	&\Eq \frac{m^2}{v^2 n}\sum_{i=1}^{\floor{n\hT}}\int_{\IR_+}
           \IE \bbbcls{ D^2f(\tX_n)\bbcls{J_{i/n,y},\IE\bcls{X_n\sp{i/n}-X_n+J_{i/n,y}}}}G(dy).\label{GGI-ADB8}
}
From Lemma~\ref{lem:gginfdiffmean} (and using the notation there), we have
\ban{
\IE\bbcls{&X_n\sp{i/n}- X_n +J_{i/n,y}}  
	\Eq J_{i/n,y}+ \frac{v^2-m^2}{m^2} I_{i/n}\oG_{i/n} + A_{n,i}.\label{GGI-ADB9}
}
\ignore{
	+ \mathbf{1}[i> n^\eps] \sum_{j=\floor{i-n^\eps}+1}^{\floor{i+n^\eps}}B_{ij}\bclr{ I_{j/n}-I_{(j+1)/n}}\\
&\qquad+J_i'\bclc{n^{-\eps(r-2)} +n^\eps\bclr{\Delta_{n^{\eps-1}}(G)+\mathbf{1}[i \leq  n^\eps]}} ,
where $\Delta_{n^{\eps-1}}(G):= \sup_{0\leq x \leq \hT-n^{\eps-1} }\bclr{G(x+n^{\eps-1})- G(x)}$,
$\norm{J_{i}'} \leq C$ and $\abs{B_{ij}} \leq Cn^\eps$. 
}
Hence \Ref{GGI-ADB7} -- \Ref{GGI-ADB9} imply that the contribution from Equation~\eq{eq:stnbdt1} is bounded by
\ban{
\bbbabs{\IE \bbbcls{ \int& D^2 f_g(\tX_n)[\hJ_{t,y}\uts] \Lambda_n(dt,dy) 
       - \int D^2 f_g(\tX_n)\bcls{J_{t,y},\IE[X_n\sp{t,y}-X_n+J_{t,y}]}\frac{\lambda_n(dt,dy)}{\sigma_n^{2}}}}\non \\
	&\Le
	\frac{C}{n}   \sum_{i=1}^{\floor{nT}} \int_{\IR_+} \bigl|\IE
                    \bcls{ D^2 f_g(\tX_n)\bcls{J_{i/n,y},A_{n,i}}}\bigr|\, G(dy) \\
                        	&{\leq C \hT
   \begin{cases}
      \norm{g}_{M'} n^{-\beta} + S_g n^{-1/2}, 
	&\text{if}\ g \mbox{ satisfies~\eqref{eq:smoothgginf}}, \\ 
	\norm{g}_{M'}\bclr{ n^{-\beta} + k^3 n^{-1}}, &\text{if}\  g(w) = F\bclr{w(t_1),\ldots, w(t_k)},
   \end{cases}}
            \label{GGI-ADB10}
}
where we have used the bounds of Lemma~\ref{lem:gginfdiffmean}, writing
\be{
(|i-\lfloor ns \rfloor| + 1)^{-(r-2)} = \sum_{j=0}^{\floor{n\hT}-1} \mathbf{1}\bcls{s\in [j/n, (j+1)/n)} (|i-j| +1)^{-(r-2)},
}
and using the same smoothness/finite number of instants arguments leading to~\eqref{eq:inttggbd1}.

\ignore{
Now, similar to bounding~\eq{eq:n2bd} above, we have that
\ba{
\babs{D^2 &f_g(w)\bcls{J_{i/n,y},I_{j/n}-I_{(j+1)/n}}} \\
	&\leq
\begin{cases}
 S_gn^{-1/2},  &g \mbox{ satisfies~\eq{eq:smoothgginf}}, \\
k^2 \norm{g}_{M'} \sum_{\ell=1}^k \mathbf{1}\bcls{t_\ell \in \big((j-1)/n, j/n\big]}, & g(w) = F\bclr{w(t_1),\ldots, w(t_k)}.
\end{cases}
}
Combining the last two displays, using~\eq{eq:2dstnbd}, and then applying the triangle inequality 
with~\eq{eq:inttggbd1} implies that~\eq{eq:stnbdt1} is bounded by
\ben{\label{GGI-1.11}
C\hT \begin{cases}
    \norm{g}_{M'}  \bclr{n^{-\eps(r-2)}+ c n^{-\beta+\eps(1+\beta)} + n^{2\eps-1}} + S_g n^{2\eps-\frac12},  
	&g \mbox{ satisfies~\eq{eq:smoothgginf}}, \\
   \norm{g}_{M'} \bclr{n^{-\eps(r-2)}+ c n^{-\beta+\eps(1+\beta)} + n^{2\eps-1}} +\norm{g}_{M'} k^3 n^{2\eps-1}, 
                   & g(w) = F\bclr{w(t_1),\ldots, w(t_k)}.
\end{cases}
}
}

To bound~\eq{eq:stnbdt2}, we first construct processes $\hX_n\st{i/n}$ and~$X_n\st{i/n}$ such that~$\hX_n\st{i/n}$
is  independent of $X_n\sp{i/n}-X_n$, the pair of processes $X_n\st{i/n}$ and~$\hX_n\st{i/n}$ are close to 
one another, and
$\law(X_n\sp{i/n}-X_n,X_n\st{i/n}) = \law(X_n\sp{i/n}-X_n,X_n)$.  To do so,
recalling the notation of Lemma~\ref{lem:gginfpalm}, let $\hMM_n\stin$ be a copy of~$M_n$ that is independent 
of both $M_n$ and~$M_n\sp{i/n}$. Define
\ba{
  \wT_{i}  &\Def \min\bbclc{n-i+1,\, \inf\bclc{j> T_i\colon M_n\bclr{(i+j)/n} = \hMM_n\stin\bclr{(i+j)/n} = 1}}, \\
  \wT_{i}' &\Def \min\bbclc{i, \inf\bclc{j> T_i'\colon M_n\bclr{(i-j)/n} = \hMM_n\stin\bclr{(i-j)/n} = 1}}.
}
Now set 
\be{
    M_n\st{i/n}(j/n) \Def
      \begin{cases}
                M_n(j/n), & i - \wT_{i}' < j <i + \wT_{i} \\
                \hMM_n\stin(j/n)  , & \mbox{otherwise}
      \end{cases}
}
and, for $(Y_j',\,j\geq1)$ i.i.d.\ with distribution function~$G$ and independent of the previous variables, set
\ba{
    N_n\st{i/n} &\Def \sum_{j= i-\wT_{i}'}^{ i+\wT_{i}} M_n\st{i/n}(j/n) \delta_{(j/n, Y_j)}
          + \sum_{j\notin {[i-\wT_{i}', i+\wT_{i}]}}  M_n\st{i/n}(j/n) \delta_{(j/n, Y_j')}; \\
    \hNN_n\stin &\Def \sum_{j = 1}^{\floor{n\hT}} \hMM_n\stin(j/n) \delta_{(j/n, Y_j')},
}
and then set
\be{
            X_n\st{i/n} \Def \int_{\cS} J_{t,y} N_n\st{i/n}(dt,dy);\qquad 
             \hXX_n\stin \Def \int_{\cS} J_{t,y} \hNN_n\stin(dt,dy).
}
It is clear that $\law(N_n\st{i/n},N_n\sp{i/n}-N_n) = \law(N_n,N_n\sp{i/n}-N_n)$, because
$N_n$ and~$N_n\st{i/n}$ differ only by having different choices of independent and identically
distributed zero--delayed renewal processes defining their continuations outside the interval
$[i-\wT_{i}',i+\wT_{i}]$, and these are independent of $N_n\sp{i/n}-N_n$, which is determined
by events defined only on $[i-\wT_{i}',i+\wT_{i}]$.

\ignore{
and also note that $N_n\st{i/n}$ is independent of $N_n\sp{i/n}-N_n$, since the latter only gives 
nonzero measure to subsets of 
\be{
\bclr{(i-T_i')/n,(i+T_i)/n}\times \IR_{+} \subseteq\bclr{(i-\wT_i')/n,(i+\wT_i)/n}\times \IR_{+},
}
on which $N_n\st{i/n}$ is determined by $\hM_n\st{i/n}$ and $(Y_j')_{j\geq1}$, which in turn are independent 
of $N_n\sp{i/n}$ and~$N_n$. Thus, $X_n\st{i/n}$ is independent of $(X_n\sp{i/n}-X_n)$.
}

Now, defining 
\be{
    \tX_n\st{i/n} \Def \sigma_n^{-1} \bclr{X_n\st{i/n}-\lambda_n}; \qquad
     \thXX_n\stin \Def \sigma_n^{-1} \bclr{\hXX_n\stin-\lambda_n},     .
}
we observe that
\ban{
    \IE \bbbcls{ \int D^2 f(\thXX_n\stin)\bcls{J_{i/n,y}, (X_n\sp{i/n}-X_n)-\IE[X_n\sp{i/n}-X_n]} G(dy)}
              \Eq 0. \label{eq:meanindeq}
}
This follows primarily because
$\hMM_n\stin$ is independent of both $M_n$ 
and~$M_n\sp{i/n}$, and hence $\thXX_n\stin$ and $X_n\sp{i/n}-X_n$ are independent.
In more detail, because
 $M_n\sp{i/n}(j/n) = M_n(j/n)$ for $j\not\in(i-T'_i,i+T'_i)$, and because of the independence of  
$\bclr{M_n\sp{i/n}-M_n}$, $(Y_j,\,j\geq1)$ and $\thXX_n\stin$, we can use Lemma~\ref{lem:intlin} to show that 
\ba{
\IE \bbbcls{& D^2 f(\thXX_n\stin)\bcls{J_{i/n,y}, (X_n\sp{i/n}-X_n)}} \\
       &=\ \IE\bbbcls{\sum_{j=i-T'_i+1}^{i+T_i-1} \bclr{M_n\sp{i/n}(j/n)-M_n(j/n)} 
                         \IE\bigl\{ D^2 f(\thXX_n\stin)\bcls{J_{i/n,y},\,I_{j/n}-I_{jn^{-1}+Y_j}}\bigr\}} \\
	&=\ \IE\bbbcls{\sum_{j=i-T'_i+1}^{i+T_i-1} \bclr{M_n\sp{i/n}(j/n)-M_n(j/n)}\int_{\IR_+} 
                          D^2 f(\thXX_n\stin)\bcls{J_{i/n,y},\,I_{j/n}-I_{jn^{-1}+y'}} G(dy')} \\
	&=\ \IE\bbbcls{\sum_{j=i-T'_i+1}^{i+T_i-1} \bclr{M_n\sp{i/n}(j/n)-M_n(j/n)}
                          D^2 f(\thXX_n\stin)\bcls{J_{i/n,y},\, I_{j/n}\overline G_{j/n}}}\\
	&=\ \sum_{j=1}^{\floor{n\hT}} \IE\bcls{M_n\sp{i/n}(j/n)-M_n(j/n)} 
               \IE\bbcls{D^2 f(\thXX_n\stin)\bcls{J_{i/n,y},\, I_{j/n}\overline G_{j/n}}}\\
	&=\ \IE\bbcls{D^2 f(\thXX_n\stin)\bbcls{J_{i/n,y},\,
                  \sum_{j=1}^{\floor{n\hT}} \IE\bcls{M_n\sp{i/n}(j/n)-M_n(j/n)} I_{j/n}\overline G_{j/n}}} \\
	&=\ \IE \bbcls{ D^2 f(\thXX_n\stin)\bcls{J_{i/n,y},\, \IE\cls{X_n\sp{i/n}-X_n}}},
}
and integrating with respect to $y$ gives~\eq{eq:meanindeq}. 
Thus, we can bound the term corresponding to~\eq{eq:stnbdt2} as follows.  
First observe that,
because $\law(N_n\st{i/n},N_n\sp{i/n}-N_n) = \law(N_n,N_n\sp{i/n}-N_n)$,
\begin{align}\label{EQ4_20}
 &\bbbabs{\frac{m^2}{v^2 n}\sum_{i=1}^{\floor{n\hT}}\int_{\IR_+}
            \IE \bbcls{D^2 f(\tX_n)\bcls{J_{i/n,y}, (X_n\sp{i/n}-X_n)-\IE[X_n\sp{i/n}-X_n]} }G(dy)} \non\\
    &\quad  =\ \bbbabs{\frac{m^2}{v^2 n}\sum_{i=1}^{\floor{n\hT}}\int_{\IR_+}
            \IE \bbcls{D^2 f(\tX_n\st{i/n})\bcls{J_{i/n,y}, (X_n\sp{i/n}-X_n)-\IE[X_n\sp{i/n}-X_n]} }G(dy)}.
\end{align}
Now we can use~\Ref{eq:meanindeq} to give
\ban{
  \eqref{EQ4_20}
    &\leq\ \frac{m^2}{v^2 n}\sum_{i=1}^{\floor{n\hT}}\int_{\IR_+} 
          \bbbabs{\IE \bbcls{D^2 f(\tX_n\st{i/n})\bcls{J_{i/n,y}, (X_n\sp{i/n}-X_n)-\IE[X_n\sp{i/n}-X_n]}  \non\\
    &\hspace{3.2cm} -D^2 f(\thXX_n\stin)\bcls{J_{i/n,y}, (X_n\sp{i/n}-X_n)-\IE[X_n\sp{i/n}-X_n]}} }G(dy) \non\\
     &\leq\ \frac{C\norm{g}_{M} }{n}\sum_{i=1}^{\floor{n\hT}}
            \IE\bbcls{\bnorm{\tX_n\st{i/n}-\thXX_n\stin}\bnorm{(X_n\sp{i/n}-X_n)-\IE[X_n\sp{i/n}-X_n]}}  \non\\
    &\leq\ \frac{C \norm{g}_{M}}{n^{3/2}}\sum_{i=1}^{\floor{n\hT}}\IE\bcls{\clr{1+ \wT_i +\wT_i'}\clr{1+ T_i + T_i'}} 
            \Le \frac{C \hT \norm{g}_M}{\sqrt{n}}, \label{GGI-1.12}
} 
where the second inequality  follows from \eq{eq:3dstnbd}, the third is because  
$N_n\st{i/n}(j/n)=\hNN_n\stin(j/n)$ for $j\notin [ i-\wT_{i}', i+ \wT_{i}]$ and
$N_n\sp{i/n}(j/n)=N_n(j/n)$ for $j\notin [ i- T_{i}', i+ T_{i}]$, and the final inequality  
is obtained  by using Cauchy-Schwarz and then noting that, by \cite[Proposition~(6.10)]{Pitman1974}, 
under the assumption~$\IE[R^3]<\infty$,  $T_i,T_i',\wT_i$ and~$\wT_i'$ all have finite
second moments, whose values depend only on~$\law(R)$.

Similarly, to bound~\eq{eq:stnbdt3}, note that
\be{
    \IE\bcls{  \norm{X_n\sp{i/n} -X_n+J_{i/n,y}}  ^2} \Le \IE\bcls{(1 + T_i + T_i')^2} \Le C,
}
\ignore{
where we again use Cauchy-Schwarz and that both of $\IE[T_i^4], \IE[(T_i')^4]$ are bounded by a constant 
depending only on~$\law(R)$; again using \cite[Proposition~(6.10)]{Pitman1974}, and the assumption 
}
again if $\IE[R^3]<\infty$. Thus
\ben{\label{GGI-1.13}
 \frac{\norm{g}_M}{2\sigma_n}\IE \bbbcls{\int_{\cS}\norm{J_{t,y}}\norm{X_n\sp{t,y}-X_n+J_{t,y}}^2
           \bclr{\sigma_n^{-2}\lambda_n(dt,dy)}} \Le \frac{C \hT \norm{g}_M}{\sqrt{n}}\,. 
}
Combining \Ref{eq:inttggbd1}, \Ref{GGI-ADB10}, \Ref{GGI-1.12} and~\Ref{GGI-1.13} yields the bound 
given in~\Ref{GGI-expec-bnd}.

\medskip
To prove the bound on the L\'evy--Prokhorov distance, we follow the template for the M/GI/$\infty$ queue, 
and use the main results of \cite{Barbour2021}, as stated in Theorem~\ref{thm:smooth} below. 
The first hypothesis of  Theorem~\ref{thm:smooth}  is satisfied, with $\kappa_1$ and~$\kappa_2$ read from the 
bound~\Ref{GGI-expec-bnd} (noting that $\norm{g}_M \leq \norm{g}_{M'}\leq \norm{g}_{M^0}$).

To bound the modulus of continuity terms, we again use Lemma~\ref{lem:modc}. To verify~\eq{eq:modch1}, for any
$0\leq s_1 < s_2 \leq \hT$, define the regions 
\ba{
 \cR_1(s_1,s_2) \Def   \cR_{s_2}\setminus \cR_{s_1}
     \quad\mbox{and}\quad \cR_2(s_1,s_2) \Def (s_1,s_2]\times \IR_+,
}
as before  {at~\eqref{eq:regions}}, so that 
\be{
  \tX_n(t)-\tX_n(s) \Eq  \tY_n(2;s,t) - \tY_n(1;s,t) ,
}
where $\tY_n(i;s_1,s_2) := \s_n^{-1/2}\{N_n(\cR_i(s_1,s_2))  - \l_n(\cR_i(s_1,s_2))\}$, $i=1,2$.
We now use Markov's inequality to bound each term in 
\be{
\IP\bclr{\abs{\tX_n(t) - \tX_n(s) }\geq \theta} 
	\Le \IP\bclr{\abs{\tY_n(1;s,t)} \geq \theta/2} 
		+\IP\bclr{\abs{\tY_n(2;s,t)} \geq \theta/2}. 
}	
First, $\s_n\tY_n(2;s,t) = \sum_{ns < i \leq nt} \bigl[ M_n(i/n) - m^{-1}\bigr]$ is the centred number of renewals
in the interval $(ns,nt]$.  By the usual renewal theory coupling arguments, as in \cite[Proposition~6.10]{Glynn1982},
writing $\hM_i := M_n(i/n)$, the sequence $(\hM_i,\,1\le i\le n\hT)$ is strong mixing  {as introduced in \cite{Rosenblatt1956}}
 with coefficients
$\alpha_j\leq k_R j^{-(r-1)}$,  {$j=1,2,\ldots,$} for a constant~$k_R < \infty$, depending only on~$\law(R)$, that we can choose to
be at least~$1$, and $\alpha_0 \Def 1/2$.  Thus, for $0 < u \le 1$, as in \cite[(1.21)]{Rio2013},
\[
  \alpha^{-1}(u) \Def \sum_{j\geq0} \bone\{ u < \alpha_j\} \Le (k_R/u)^{1/(r-1)}  {+ \bone\{u< 1/2\}} \Le 2(k_R/u)^{1/(r-1)}.
\]
Applying \cite[Theorem~2.2]{Rio2013}, it follows that, for any $l \in \IN$, there is a constant~$C_{l,R}$ depending 
only on $\law(R)$ and~$l$ such that
\ba{
    \IE&\bbcls{\bclr{\tY_n(2;s,t)}^{2l}}  \\
    &     \Le C_{l,R}\s_n^{-2l} \bbbclc{\bbbclr{ \sum_{ns < i \leq nt} \int_0^1 \alpha^{-1}(u) Q_i^2(u)\, du}^l  
   + \sum_{ns < i \leq nt}\int_0^1 \bcls{\alpha^{-1}(u)}^{2l-1}Q_i^{2l}(u)\, du },
}
where $Q_i := q_{1/m}$ and, {for $w\in [0,1]$ and $u\in(0,1]$, 
\ben{\label{quantilefunc}
   q_w(u) =\begin{cases}
   w\bone_{(0,1-w]}(u) + (1-w)\bone_{(1-w,1]}(u), & w\geq 1/2,  \\
 (1- w)\bone_{(0,w]}(u) + w\bone_{(w,1]}(u), & w< 1/2.  \\
 \end{cases} 
} 
}
Straightforward computing now shows that, for $l < r/2$, 
\be{
\IE\bbcls{\bclr{\tY_n(2;s,t)}^{2l}} \Le K_{l,R}\s_n^{-2l} {\bcls{(\lfloor nt \rfloor - \lfloor ns \rfloor)^l + (\lfloor nt \rfloor - \lfloor ns \rfloor)}},
}
for a constant $K_{l,R} < \infty$.  Markov's inequality and~\Ref{ADB-variance-def} thus imply that, 
if $n(t-s) \ge 1/2$ and for $l < r/2$,
\ben{\label{ADB8}
   \IP\bclr{\abs{\tY_n(2;s,t) }\geq \sigma_{n} \theta/2} 
                  \Le K_{l,R}\Bigl(\frac{m^3}{v^2}\Bigr)^l 2^{2l {+1}}\theta^{-2l}3^l(t-s)^{l} \ =:\ C_l\ui\theta^{-2l}(t-s)^l.
}

 {
For $\tY_n(1;s,t)$, we observe that
\be{
     {N}_n(\cR_1(s,t))-\lambda_n(\cR_1(s,t)) \Eq \sum_{1\leq i  \leq nt} \big[ M_n(i/n) B_{i,n} - m^{-1} p_{i,n} \big], 
}
where 
\be{
   B_{i,n} \Def \mathbf{1}[(s- in^{-1})_+ <  Y_i \leq t - in^{-1}  ]\ \sim\ \text{Be}(p_{i,n});\quad   
     p_{i,n} \Def G(t - i/n) - G((s - i/n)_+).
}
The fact that the random variables $(B_{i,n},\,1\le i\le n\hT)$ are independent of~$N_n$ implies that the mixing properties
of the sequence $(\hM_i,\,1\le i\le n\hT)$ are inherited by the sequence $(\hM_i B_{i,n},\,1\le i\le n\hT)$,
so that \cite[Theorem~2.2]{Rio2013} can be applied with the same function~$\a^{-1}$, giving 
\besn{
\IE\bbcls{&\bclr{\tY_n(1;s,t)}^{2l}} \\
	&\Le C_{l,R}\s_n^{-2l}
           \bbbclc{\bbbclr{ \sum_{1\leq i\leq nt} \int_0^1 \alpha^{-1}(u) Q_i^2(u) du}^l  
    + \sum_{1\leq i \leq nt}\int_0^1 \bcls{\alpha^{-1}(u) }^{2l-1}Q_i^{2l}(u) du }, \label{ADB7a}
}
where now $Q_i := q_{p_{i,n}/m}$.  Using~\eqref{quantilefunc}, we have that for $1 \le l < r/2$,
\besn{
\int_0^1 (\alpha^{-1}(u))^{2l-1} q_w^{2l}(u)\, du
    &           \Le c_{1}(r,R,l) w^{(r-2l)/(r-1)}. \label{ADB8a}
}
From the assumptions on~$G$, by comparing sums and integrals, it follows that, for $1/(2n) \le (t-s) \le 1$,
\ban{
     n^{-1}\sum_{1 \le i \le nt} p_{i,n}  &\Le \int_s^t G(v)\,dv + 2n^{-1} \Le
          5(t-s), \label{ADB9}
}
and, using~\Ref{ADB8a}, for $1 \le l < r/2$ such that $(r-2l)/(r-1) \ge \LL$, we have
\ban{
     n^{-1}\sum_{1 \le i \le nt} p_{i,n}^{(r-2l)/(r-1)}  &\Le (t-s)^{\b(r-2l)/(r-1)}
            \Bigl\{\int_0^\infty \bigl(g_G(v)\bigr)^{(r-2l)/(r-1)}\,dv + 2\bigl(g_G(0)\bigr)^{(r-2l)/(r-1)} \Bigr\}\non\\
             &\Le c_3(r,R,l)(t-s)^{\b(r-2l)/(r-1)} .     \label{ADB9a}
}
Hence, from~\Ref{ADB7a}, if $r(1-\b) \ge 1$, writing $\b_r := \b(r-2)/(r-1)$, we find that,
for $0 \le s < t \le \hT$ and for $n(t-s) \ge 1/2$, 
\be{
  \IE\bbcls{\bclr{\tY_n(1;s,t)}^{2l}} 
      \Le C_{l,R}(n/\s_n^2)^{l} (t-s)^{l\b_r},
}
for any $l < (r - \LL(r-1))/2$.  {The assumption $r(1-\b) \ge 1$ ensures the exponent $l\b_r$ is no larger than those appearing when applying~\eqref{ADB7a},~\eqref{ADB8a} and~\eqref{ADB9a}.} It now follows, by Markov's inequality and~\Ref{ADB-variance-def}, that,
for such~$l$, and if $r(1-\b) \ge 1$,
\ben{\label{ADB10}
    \IP\bclr{\abs{\tY_n(1,s,t)}\geq  \theta/2} \Le C_l\ui  \theta^{-2l}(t-s)^{l\b_r},\quad n(t-s) \ge 1/2.
}
If $r(1-\b) < 1$, the inequality~\Ref{ADB10} holds only for $0 \le s < t \le \hT$ such that $(t-s)n^{(r-1)/r\b} \ge 1/2$.
}

To verify~\eq{eq:modch2} for $\tX_n$, we note that, for $(k-1)/\MM \le u \le k/\MM$,
\ba{
   |\tY_n(i;(k-1)/\MM,u)| &\Le |\tY_n(i,(k-1)/\MM,k/\MM)|
       + 2\s_n^{-1}\l_n(\cR_i((k-1)/\MM,k/\MM)),
}
$i=1,2$, where $\l_n(\cR_2((k-1)/\MM,k/\MM)) = m^{-1}n/\MM$ and, because of~\Ref{ADB9}, 
$\l_n(\cR_1((k-1)/\MM,k/\MM)) \le 5m^{-1}n/\MM$.  Writing $\rho := \rho(r,\b) := \min\{1,(r-1)/r\b\}$, and
taking $\MM := n^{\rho}$, it follows
from \Ref{ADB8} and~\Ref{ADB10} that  
 we can take
\ben{\label{ADB11}
               \f_{\MM}(\th) \Def {\MM}^{1-l\b_r}C_l\uh\th^{-2l}
}
in~\eq{eq:modch2}, for a suitable constant~$C_l\uh$, if $\th > 48m^{-1}n^{1-\rho}/\s_n$. 
We can now apply Lemma~\ref{lem:modc} with $b=2l$ and $a = l\b_r-1$, for 
$l = l_r := \lceil (r - \LL(r-1))/2 \rceil - 1$, and with $\MM := n^{\rho}$,
 to find that, for $\theta > 48\bigl(1-2^{-(l_r\b_r-1)/4l_r}\bigr)^{-1} m^{-1}n^{1-\rho}/\sigma_n$ 
and $\epsilon\in (n^{-\rho}, 1)$, 
\ban{
  \IP(\omega_{\tX_n} (\epsilon) \geq \theta/2) 
             \Le    C\hT  \th^{-2l_r} \ee^{l\b_r - 1}, \label{modctxn}
}
for some suitable constant~$C$.

For the modulus of continuity of~$Z$, an easy calculation shows that for any $0\leq r < s \leq \hT$, 
\be{
  \IE\bcls{\bclr{Z(s)-Z(r)}^2}\leq  C (s-r)^\beta.
}
{and so \cite[Remark~1.6]{Barbour2021}
implies  that for any $l \ge 1$,
\ben{\label{modczgg}
   \IP[\omega_{Z} (\epsilon) \geq \theta]  \Le   C\hT \theta^{-{2l}} \epsilon^{l\bo-1}.
   }
}

 {
From \Ref{eq:inttggbd1}, \Ref{GGI-ADB10}, \Ref{GGI-1.12} and~\Ref{GGI-1.13}, and for  $\hT \ge 1$,
we can now apply Theorem~\ref{thm:smooth} with
\bes{
  \kappa_1
	& =\  \bigo\bclr{\hT n^{-\bab}}, \quad\mbox{where}\quad 
                 \bab \Def \min\{\beta,1/2\};\qquad
  \kappa_2 \Eq \bigo\bclr{\hT n^{-1/2}}.
}
This, using \Ref{modctxn} and~\Ref{modczgg}, for any choice of $\ee,\d > 0$ and $\th > cn^{1/2 - \rho(r,\b)}$, 
and for any $l \ge 1$, implies a bound of 
\[
 C    \bigl( \d\sqrt{ {\hT} \log n} + \th + \hT^{5/2} n^{-\bab}(\ee\d)^{-3} + \hT^{ {3/2}} n^{-1/2}(\ee\d)^{-2}
        + \hT \ee^{l_r\b_r - 1}\th^{-2l_r}  + \hT \ee^{l \b-1}\th^{-2l} \bigr),
\]
for a suitable constant~$C$, where we recall that $\b_r := \b(r-2)/(r-1)$ and that 
$\rho(r,\b) := \min\{1,(r-1)/r\b\}$.  Taking
\ba{
    \th &\Eq  {\sqrt{\hT}}\d \Eq \bigl\{ (\hT^{4} n^{-\bab})^{l_r\b_r-1} \hT^{ {3}} \bigr\}^{1/(6l_r + 4l_r\b_r -1)}; \quad
       \ee  \Eq \bigl\{ (\hT^{ {4}} n^{-\bab})^{2l_r+1} \hT^{ {-4}} \bigr\}^{1/(6l_r + 4l_r\b_r -1)},
}
gives a bound of order
\[
     \bigo\Bigl( \sqrt{\log n}\,\Bigl\{ (\hT^{5/2} n^{-\bab})^{l_r\b_r-1} \hT^{ {3}} \Bigr\}^{1/(6l_r + 4l_r\b_r -1)} \Bigr).
\] 
Here, we note that  {if $\epsilon\leq 1$, for example if $\hT\leq n^{\psi'}$ for $\psi'< \bar\beta(2l_r+1)/(8l_r)$, then the term $\hT \ee^{l\b-1}\th^{-2l}$ can be made smaller order than the others by choosing $l=l_r$, noting $\beta_r < \beta$.
}
A calculation also shows that this choice of~$\th$ indeed satisfies
$\th > 48(1-2^{-(l_r\b_r - 1)/4l_r})^{-1} m^{-1}n^{1-\rho}/\sigma_n$ for all~$n$ sufficiently large. 
}
\end{proof}

\section{Smoothing result}

We state a specific consequence of the main results of \cite{Barbour2021}, which we use above to prove weak convergence.
Let $M^0\subset M'$ be the set of functions $h\colon \ID{^p}\to \IR$ such that
\be{
 \norm{h}_{M^0}\Def \sup_{w\in \ID^p} \abs{h(w)} +\sup_{w\in \ID^p}\norm{D h(w)} +\sup_{w\in \ID^p}\norm{D^2 h(w)}
      +\sup_{w,v\in \ID^p}\frac{\norm{D^2 h(w+v)-D^2 h(w)}}{\norm{v}}
}
is finite. Note that $\norm{h}_{M'}\leq \norm{h}_{M^0}$.
For $x\in\ID^p$, let $\omega_x(\eps):= \sup_{0\leq s,t\leq \hT: \abs{s-t} < \eps} \abs{x(t)-x(s)}$ denote the modulus of 
continuity of~$x$.
 
\begin{theorem}[Theorem~1.1 of \cite{Barbour2021}]\label{thm:smooth}
Let $Y, Z$ be random elements of {$\ID^p := \ID([0,\hT], \IR^p)$}, with $\hT\geq 1$, such that~$Z$ has almost surely 
continuous sample paths. 
Suppose that there are $\kappa_1,\kappa_2\geq 0$ such that for  any $g\in M^0$ satisfying the smoothness 
condition~\eq{eq:smoothgginf},  we have
\be{
 |\IE g(Y) - \IE g(Z)| \ \leq\  \kappa_1 \norm{h}_{M^0} + S_g \, \kappa_2.
}
Letting 
$\dlp$ denote the L\'evy-Prokhorov metric, 
we have that for any positive 
$\delta, \epsilon, \theta, \gamma$ with $\epsilon,\delta\in (0,1)$, 
\ba{
  \dlp\bclr{\law(&X), \law(Z)} \\
      &\leq\ \hC  \max\bbclc{\theta+\gamma,  
      \frac{\kappa_1\hT^{3/2}}{(\epsilon \delta)^{3}} +\frac{\kappa_2 \hT^{1/2}}{ (\epsilon \delta)^{2} }
           + \IP\bclr{\omega_Y(\epsilon)\geq \theta}  + \IP\bclr{\omega_Z(\epsilon)\geq \theta} +
	     e^{-\frac{\gamma^2}{2{p}\hT\delta^2}} },
}
where $\hC$ is a universal constant.
\end{theorem}

To bound the modulus of continuity terms appearing in the previous theorem, we use the following lemma,
also noted in \cite[Lemma~1.3 and Remark 1.4(1)]{Barbour2021}, {applied to each component.}

\begin{lemma}\label{lem:modc}
Let $X\in \ID$ be such that there are positive constants $a, b$ and $K$ such that
\ben{\label{eq:modch1}
   \IP\bclr{|X(s) - X(t)| \geq \theta} \Le K \theta^{-b}|s-t|^{1+a} \ \mbox{ for }\ 
                                     \tsfrac{1}{2} \MM^{-1} \leq |s-t| \leq 1/2,
}
and that 
\ben{\label{eq:modch2}
 \MM \, \IP\bbclr{\sup_{(k-1)/\MM \le s \le k/\MM}|X(s) - X\bclr{(k-1)/\MM}| \geq \theta} \Le \phi_\MM(\theta)\ 
         \mbox{ for }\ 1 \le k \le \lceil \MM \hT \rceil.  
}
Then, for any $\epsilon \in (\MM^{-1},1)$, 
\[
     \IP\bclr{\omega_{X}(\epsilon) \ \geq\ \theta}
             \Le 2 \hT\bbclc{\phi_\MM\bclr{\theta\bigl(1-2^{-a/(2b)}\bigr)/  18 } + C'( K,a, b) \theta^{-b} \epsilon^{a}},
\]
for a constant $C'(K,a, b) < \infty$.
\end{lemma}

We have also used the following, technical lemma.

\begin{lemma}\label{lem:intlin} 
Assume $f\in M'$ either satisfies the smoothness condition~\eq{eq:smoothgginf}, or is a function of the form 
$f(w)=F\bclr{w(t_1),\ldots, w(t_k)}$ for some {$F\colon(\IR^p)^k \to \IR$} and $\{t_1,\ldots,t_k\}\subseteq[0,\hT]$. 
Letting $I_r(s):=\bone\{s\geq r\}$ and $G_t(s):=G(s-t)$ for a distribution function~$G$, 
then {for any $w,x_1,x_2\in \ID^p$ and $r,t\geq0$, we have 
\ba{
\int_{\IR_+} D^2 f(w)[x_1I_r, x_2I_{t+y}] G(dy)
	&\Eq D^2 f(w)\bbcls{x_1I_r, \int_{\IR_+} x_2I_{t+y}G(dy)}  \Eq  D^2 f(w)\bcls{x_1I_r, x_2G_t }.
}
and
\ba{
\int_{\IR_+} D^2 f(w)[x_1G_r, x_2I_{t+y}] G(dy)
	&\Eq D^2 f(w)\bbcls{x_1G_r, \int_{\IR_+} x_2I_{t+y}G(dy)} \\
	&  \Eq  D^2 f(w)\bcls{x_1G_r, x_2G_t }.
}
}
\end{lemma}
\begin{proof}
If $f(w)=F\bclr{w(t_1),\ldots, w(t_k)}$, then a simple calculation shows that 
\be{
   D^2 f(w)[w_1,w_2] \Eq \sum_{i,j=1}^k {w_1(t_i)^\tr F_{ij}\bclr{w(t_1),\ldots, w(t_k)}  w_2(t_j)}, 
}
where we write $F_{ij}$ for the { $p\times p$ matrix corresponding to the } mixed partial of $F$ in coordinates $i$ and $j$. The result now follows directly 
after noting that $\int I_{t+y} G(dy) = G_t$.

Now assume that~$f$ satisfies the smoothness condition~\eq{eq:smoothgginf}.
Both results are obviously true from bilinearity if $G$ is a discrete distribution function. By considering the atoms 
separately, we can without loss of generality assume $G$ is continuous. We show the result by approximating $G$ by a 
discretised version. Let $Y\sim G$ and define $Y_m:=\floor{m Y}/m$,  noting that 
$Y_m$ converges almost surely (so in distribution) to $Y$, as $m\to \infty$.
For the first assertion, the function
{
\be{
   y\ \mapsto\ D^2 f(w)[x_1I_r, x_2I_{t+y}] 
}
}is continuous for $y\geq0$, since
{
\ba{
  D^2 f(w)[x_1I_r, x_2I_{t+y+\eps}]-D^2 f(w)[x_1I_r, x_2I_{t+y}]
	&\Eq D^2 f(w)\cls{x_1I_r, x_2(I_{t+y+\eps}-I_{t+y})} \to 0,
}
by~\eq{eq:smoothgginf}.
}
By applying Lemma~\ref{lem:bilnbds} and noting that $f\in M'$, it is also easy to see the function is bounded. Thus, 
because of weak convergence, as $m\to\infty$,
{
\be{
  \IE\bcls{D^2 f(w)[x_1I_r, x_2I_{t+Y_m}]}\ \to\ \IE\bcls{D^2 f(w)[x_1I_r, x_2I_{t+Y}]}
               \Eq \int_{\IR_+} D^2 f(w)[x_1I_r, x_2I_{t+y}] G(dy).
}
}
On the other hand, because $Y_m$ is discrete, we have
{
\be{
   \IE\bcls{D^2 f(w)[x_1I_r, x_2I_{t+Y_m}]} \Eq D^2 f(w)\bcls{x_1I_r, x_2\IE\cls{I_{t+Y_m}}},
}
}
and, again because of Lemma~\ref{lem:bilnbds} and because $f\in M'$, we have
{
\ba{
   \bbabs{D^2 f(w)\bcls{x_1I_r, x_2\IE\cls{I_{t+Y_m}}}-D^2 f(w)&\bcls{x_1I_r, x_2\IE\cls{I_{t+Y}}}}  \\
                       &  \Le 3 \| f\|_{M'}\|x_1\|\,\|x_2\| \bnorm{\IE\cls{I_{t+Y_m}}-\IE\cls{I_{t+Y}}}
}
which converges to zero,
}
since
\ba{
\bnorm{\IE\cls{I_{t+Y_m}}-\IE\cls{I_{t+Y}}}
	&\Eq \sup_{y\in[0,\hT]} \babs{\IP(Y\leq y-t)-\IP(Y_m\leq y-t)} \\
	&\Le \sup_{y\in[0,\hT]} \bbabs{\IP\bbclr{Y\in \bbclr{y, \frac{\floor{my} +1}{m}}}} \to 0,
                ~\text{by the continuity of $G$.}
}

For the second assertion, the function {$y \longmapsto D^2 f(w)[x_1G_r, x_2I_{t+y}]$}
is bounded continuous on $\mathbb{R}_+$, since, using the first assertion as well as  the   
condition~\eq{eq:smoothgginf},
{
\ba{
&D^2 f(w)[x_1G_r, x_2I_{t+y+\eps}]-D^2 f(w)[x_1G_r, x_2I_{t+y}] \\
	&\Eq \int_{\IR_+} D^2 f(w)\cls{x_1I_{r+y'}, x_2(I_{t+y+\eps}-I_{t+y})}G(dy') \xrightarrow{\eps\downarrow 0   } 0.
}
}The rest of the proof follows in exactly the same way as for the first assertion, replacing $I_r$ by~$G_r$.
\end{proof}
%
%

%
%

\end{document}